\documentclass{article}

\usepackage{amsmath,amssymb}

\DeclareMathOperator{\coker}{coker}

\DeclareMathOperator{\fchar}{char}

\DeclareMathOperator{\Gal}{Gal}
\DeclareMathOperator{\Hom}{Hom}

\DeclareMathOperator{\ord}{ord}
\DeclareMathOperator{\Pic}{Pic}
\DeclareMathOperator{\rk}{rk}
\DeclareMathOperator{\Spec}{Spec}

\DeclareMathOperator{\vol}{vol}

\newcommand{\CC}{\mathbb{C}}
\newcommand{\FF}{\mathbb{F}}

\newcommand{\QQ}{\mathbb{Q}}
\newcommand{\RR}{\mathbb{R}}
\newcommand{\ZZ}{\mathbb{Z}}

\renewcommand{\AA}{\mathbb{A}}

\renewcommand{\det}{\operatorname{det}}

\renewcommand{\Re}{\operatorname{Re}}
\renewcommand{\emptyset}{\varnothing}

\newcommand{\ar}{\text{\it ar}}
\newcommand{\BM}{\text{\it BM}}

\newcommand{\et}{\text{\it \'{e}t}}
\newcommand{\fg}{\text{\it fg}}

\newcommand{\red}{\text{\it red}}
\newcommand{\tors}{\text{\it tors}}
\newcommand{\Wc}{\text{\it W,c}}
\newcommand{\Zar}{\text{\it Zar}}

\newcommand{\dfn}{\mathrel{\mathop:}=}

\newcommand{\RHom}{R\!\Hom}

\usepackage{tikz-cd}
\usetikzlibrary{arrows}
\usetikzlibrary{calc}
\usetikzlibrary{babel}
\usetikzlibrary{decorations.pathreplacing}
\usetikzlibrary{patterns}

\usepackage{hyperref}

\hypersetup{
  final,
  colorlinks,
  linkcolor={red!60!black},
  citecolor={blue!60!black},
  urlcolor={blue!80!black}
}

\usepackage{amsthm}

\newtheoremstyle{myplain}
{\topsep}   
{\topsep}   
{\itshape}  
{0pt}       
{\scshape} 
{.}         
{5pt plus 1pt minus 1pt} 
{}   

\theoremstyle{myplain}

\newtheorem*{maintheorem*}{Main theorem}
\newtheorem*{thetheorem*}{Theorem}
\newtheorem*{proposition*}{Proposition}
\newtheorem{theorem}{Theorem}[section]
\newtheorem{proposition}[theorem]{Proposition}
\newtheorem{lemma}[theorem]{Lemma}

\newtheoremstyle{mydefinition}
{\topsep}   
{\topsep}   
{}  
{0pt}       
{\scshape} 
{.}         
{5pt plus 1pt minus 1pt} 
{}   

\theoremstyle{mydefinition}
\newtheorem*{conjecture*}{Conjecture}
\newtheorem{definition}[theorem]{Definition}

\newtheorem{remark}[theorem]{Remark}
\newtheorem{example}[theorem]{Example}

\title{Zeta-values of one-dimensional arithmetic schemes at strictly negative integers}
\author{Alexey Beshenov}
\date{}

\AtEndDocument{%
  \par
  \medskip
  \begin{tabular}{@{}l@{}}%
    \\
    Alexey Beshenov \\
    E-mail: \texttt{cadadr@gmail.com} \\
    URL: \url{https://cadadr.org/}
  \end{tabular}}

\begin{document}

\maketitle

\begin{abstract}
  Let $X$ be an arithmetic scheme (i.e., separated, of finite type over
  $\Spec \ZZ$) of Krull dimension $1$. For the associated zeta function
  $\zeta (X,s)$, we write down a formula for the special value at $s = n < 0$ in
  terms of the \'{e}tale motivic cohomology of $X$ and a regulator. We prove it
  in the case when for each generic point $\eta \in X$ with
  $\fchar \kappa (\eta) = 0$, the extension $\kappa (\eta)/\QQ$ is
  abelian. We conjecture that the formula holds for any one-dimensional
  arithmetic scheme.

  This is a consequence of the Weil-\'{e}tale formalism developed by the author
  in \cite{Beshenov-Weil-etale-1} and \cite{Beshenov-Weil-etale-2}, following
  the work of Flach and Morin \cite{Flach-Morin-2018}. We also calculate the
  Weil-\'{e}tale cohomology of one-dimensional arithmetic schemes and show that
  our special value formula is a particular case of the main conjecture from
  \cite{Beshenov-Weil-etale-2}.
\end{abstract}



\section{Introduction}

Let $X$ be an \textbf{arithmetic scheme}, by which we mean in this text that it
is separated and of finite type over $\Spec \ZZ$.  The \textbf{zeta function}
associated to $X$ (see, e.g. \cite{Serre-1965}) is given by
\[ \zeta (X,s) \dfn \prod_{\substack{x \in X \\ \text{closed pt.}}}
  \frac{1}{1 - N (x)^{-s}}, \]
where the norm of a closed point $x\in X$ is the size of the corresponding
residue field:
$$N (x) \dfn |\kappa (x)| \dfn |\mathcal{O}_{X,x}/\mathfrak{m}_{X,x}|$$
The above product converges for $\Re s > \dim X$ and is supposed to have a
meromorphic continuation to the whole complex plane. Although the latter is a
wide-open conjecture in general, it is well-known for one-dimensional schemes,
which is the case of interest in this article.

If $\zeta (X,s)$ admits a meromorphic continuation around $s = n$, we denote by
\begin{equation}
  \label{eqn:vanishing-order}
  d_n \dfn \ord_{s=n} \zeta (X,s)
\end{equation}
the \textbf{vanishing order} of $\zeta (X,s)$ at $s = n$. The corresponding
\textbf{special value} of $\zeta (X,s)$ at $s = n$ is defined as the leading
nonzero coefficient of the Taylor expansion:
$$\zeta^* (X,n) \dfn \lim_{s \to n} (s - n)^{-d_n}\,\zeta (X,s).$$

Since the 19th century, many formulas (both conjectural and unconditional) have
been proposed to interpret the numbers $\zeta^* (X,n)$ in terms of geometric and
algebraic invariants attached to $X$. A primordial example is Dirichlet's
\textbf{analytic class number formula}. For a number field $F/\QQ$, we denote by
$\mathcal{O}_F$ the corresponding ring of integers. Then
$$\zeta_F (s) \dfn \zeta (\Spec \mathcal{O}_F, s)$$
is the \textbf{Dedekind zeta function} attached to $F$. From the well-known
functional equation for $\zeta_F (s)$, it is easy to see that it has a zero at
$s = 0$ of order $r_1 + r_2 - 1$, where $r_1$ (resp.  $2 r_2$) is the number of
real embeddings $F \hookrightarrow \RR$ (resp. complex embeddings
$F \hookrightarrow \CC$). The corresponding special value at $s = 0$ is given by
\begin{equation}
  \label{eqn:zeta-F-at-s=0}
  \zeta^*_F (0) = -\frac{h_F}{\omega_F}\,R_F,
\end{equation}
where $h_F = |\Pic (\mathcal{O}_F)|$ is the class number,
$\omega_F = |(\mathcal{O}_F)^\times_\tors|$ is the number of roots of unity in
$F$, and $R_F \in \RR$ is the regulator. See, e.g.,
\cite[Chapter~5, \S 1]{Borevich-Shafarevich} or \cite[\S VII.5]{Neukirch-1999}.

The question naturally arises whether there are formulas similar to
\eqref{eqn:zeta-F-at-s=0} for $s = n \in \ZZ$ other than $s = 0$ (or $s = 1$,
which is related to $s = 0$ via the functional equation). To do this, one must
find a suitable generalization for the numbers $h_F$, $\omega_F$, $R_F$.  Many
special value conjectures of varying generality go back to this question.

Lichtenbaum proposed formulas in terms of algebraic $K$-theory in his pioneering
work \cite{Lichtenbaum-1973}. Later these were also reformulated in terms of
$p$-adic cohomology $H^i (\Spec \mathcal{O}_F [1/p]_\et, \ZZ_p (n))$ for
$i = 1,2$ and all primes $p$; the corresponding formula is known as the
\textbf{cohomological Lichtenbaum conjecture}; see, for example,
\cite[\S 1.7]{Huber-Kings-2003} for the statement and a proof for abelian number
fields $F/\QQ$. We will not go into details here, since it is more convenient
for us to use motivic cohomology instead of working with $p$-adic cohomology for
varying $p$.

A suitable generalization of $R_F$ are the \textbf{higher regulators} considered
since the work of Borel \cite{Borel-1977} and later by Beilinson
\cite{Beilinson-1984}.

We do not attempt to give an adequate historical survey of the subject or to
write down all the conjectured formulas; the interested reader may consult,
e.g., \cite{Kolster-2004,Goncharov-2005,Kahn-2005}.

\vspace{1em}

Later, Lichtenbaum proposed a new research program known as
\textbf{Weil-\'{e}tale cohomology}; see
\cite{Lichtenbaum-2005,Lichtenbaum-2009-Euler-char,Lichtenbaum-2009-number-rings,Lichtenbaum-2021}.
It suggests that for an arithmetic scheme $X$ the special value of $\zeta (X,s)$
at $s = n \in \ZZ$ can be expressed in terms of the Weil-\'{e}tale cohomology,
which is a suitable modification of the \'{e}tale motivic cohomology of $X$.
Flach and Morin in \cite{Flach-Morin-2018} gave a construction of Weil-\'{e}tale
cohomology groups $H^i_\Wc (X,\ZZ(n))$ for a proper and regular arithmetic
scheme $X$, and stated a precise conjectural relation of $H^i_\Wc (X,\ZZ(n))$ to
the special value $\zeta^* (X,n)$.

In \cite[\S 5.8.3]{Flach-Morin-2018} they write down an explicit formula for the
case of $X = \Spec \mathcal{O}_F$. For $n \le 0$ and in terms of cohomology
groups $H^i (X_\et, \ZZ^c (n))$, it reads
\begin{equation}
  \label{eqn:flach-morin-zeta-F-formula}
  \zeta_F^* (n) = \pm\frac{|H^0 (X_\et, \ZZ^c (n))|}{|H^{-1} (X_\et, \ZZ^c (n))_\tors|}\,R_{F,n}
  \quad \text{for }n \le 0.
\end{equation}
The definition of $H^i (X_\et, \ZZ^c (n))$ is reviewed below.
The regulator $R_{F,n} = R_{\Spec \mathcal{O}_F,n}$ is defined in
\S\ref{sec:regulators}.

By \cite[Proposition~5.35]{Flach-Morin-2018}, formula
\eqref{eqn:flach-morin-zeta-F-formula} holds unconditionally for abelian number
fields $F/\QQ$, via a reduction to the \textbf{Tamagawa number conjecture} of
Bloch--Kato--Fontaine--Perrin-Riou.

In particular, if we take $n = 0$, then $\ZZ^c (0) \cong \mathbb{G}_m [1]$, and
$R_{F,0}$ is the usual Dirichlet regulator, so
\eqref{eqn:flach-morin-zeta-F-formula} becomes the classical formula
\eqref{eqn:zeta-F-at-s=0}:
\[ \zeta_F^* (0) =
  \pm \frac{|H^1 (\Spec \mathcal{O}_{F,\et}, \mathbb{G}_m)|}{|H^0 (\Spec \mathcal{O}_{F,\et}, \mathbb{G}_m)_\tors|}\,R_F =
  \pm \frac{|\Pic (\mathcal{O}_F)|}{|(\mathcal{O}_F)^\times_\tors|}\,R_F, \]

We also mention that Flach and Morin have a similar special value formula for
$n > 0$, which includes a correction factor $C (X,n) \in \QQ$. In this text we
will say nothing about the case of $n > 0$; the reader can consult
\cite{Flach-Morin-2018} for more details, and also the subsequent papers
\cite{Flach-Morin-2020,Flach-Morin-2020-Muenster,Morin-2021-THH} which shed
light on the nature of the correction factor $C (X,n)$.

For $n < 0$, the author in \cite{Beshenov-Weil-etale-1} and
\cite{Beshenov-Weil-etale-2} extended the work of Flach and Morin
\cite{Flach-Morin-2018} to an arbitrary arithmetic scheme $X$ (thus removing the
assumption that $X$ is proper or regular). In this text, we would like to work
out explicitly the corresponding special value formula for one-dimensional
arithmetic schemes.

\vspace{1em}

To state the main result, it is useful to introduce the following terminology.

\begin{definition}
  \label{dfn:abelian-scheme}
  We say that a one-dimensional arithmetic scheme $X$ is \textbf{abelian} if
  each generic point $\eta \in X$ with $\fchar \kappa (\eta) = 0$ corresponds to
  an abelian extension $\kappa (\eta)/\QQ$.
\end{definition}

If $X$ lives in positive characteristic, then it is trivially abelian.
The term ``abelian'' is ad hoc and was suggested by analogy with the notion of
\textbf{abelian number fields}. Hopefully there is no confusion with the
``abelian schemes'' that are generalizations of abelian varieties.

\pagebreak

Our goal is to prove the following result.

\begin{theorem}
  \label{main-theorem}
  For an abelian one-dimensional arithmetic scheme $X$, the special value of
  $\zeta (X,s)$ at $s = n < 0$ is given by
  \begin{equation}
    \label{eqn:special-value-formula}
    \zeta^* (X,n) =
    \pm 2^\delta\,\frac{|H^0 (X_\et, \ZZ^c (n))|}{|H^{-1} (X_\et, \ZZ^c (n))_\tors| \cdot |H^1 (X_\et, \ZZ^c (n))|}\,R_{X,n}.
  \end{equation}

  Here
  \begin{itemize}
  \item $H^i (X_\et, \ZZ^c (n))$ the \'{e}tale motivic cohomology
    from \cite{Geisser-2010};

  \item the correction factor $2^\delta$ is given by
    \begin{equation}
      \label{eqn:delta}
      \delta = \delta_{X,n} =
      \begin{cases}
        r_1, & n \text{ even}, \\
        0, & n \text{ odd},
      \end{cases}
    \end{equation}
    where $r_1 = |X (\RR)|$ is the number of real places of $X$,

  \item $R_{X,n}$ is a positive real number defined
    via a regulator map in \S\ref{sec:regulators}.
  \end{itemize}
\end{theorem}

We further conjecture that formula \eqref{eqn:special-value-formula} holds for
all one-dimensional arithmetic schemes, not necessarily abelian. This is
equivalent to the Tamagawa number conjecture for non-abelian number fields
(see Remark~\ref{rmk:TNC}).

\vspace{1em}

We give two proofs of \eqref{eqn:special-value-formula}: first a direct argument
in \S\ref{sec:direct-proof} and then an argument in terms of Weil-\'{e}tale
cohomology in \S\ref{sec:Weil-etale-proof}. In fact, we note that the special
value formula is the same as the conjecture $\mathbf{C} (X,n)$ formulated in
\cite{Beshenov-Weil-etale-2}, which is specialized to one-dimensional $X$ and
spelled out explicitly.

The purpose of this text is twofold. First, we establish a new special value
formula, which generalizes several formulas found in the literature. Second, we
review the construction of Weil-\'{e}tale cohomology $H^i_\Wc (X, \ZZ(n))$ from
\cite{Beshenov-Weil-etale-1} and the special value conjecture from
\cite{Beshenov-Weil-etale-2} and explain it in the case of one-dimensional
schemes. It is not very surprising that a special value formula like
\eqref{eqn:special-value-formula} exists, but the right cohomological invariants
to state it have been suggested by the Weil-\'{e}tale framework.

This text was inspired in part by the work of Jordan and Poonen
\cite{Jordan-Poonen-2020}, which deals with a formula for $\zeta^* (X,1)$, where
$X$ is an affine reduced one-dimensional arithmetic scheme. The affine and
reduced constraint does not appear in our case because work with different
invariants. Since $\zeta (X,s) = \zeta (X_\red, s)$, the ``right'' invariants
should not distinguish between $X$ and $X_\red$, and motivic cohomology
satisfies this property.

\subsection*{Notation and conventions}

\paragraph{Abelian groups.}
For an abelian group $A$, we denote
\begin{align*}
  A^D & \dfn \Hom (A, \QQ/\ZZ), \\
  A^* & \dfn \Hom (A, \ZZ).
\end{align*}
There is an exact sequence
\begin{equation}
  \label{eqn:exact-sequence-A*-etc}
  0 \to A^* \to \Hom (A,\QQ) \to A^D \to (A_\tors)^D \to 0
\end{equation}

Note that for a finite rank group $A$, the $\ZZ$-dual $A^*$ is free and has the same
rank. If $A$ is finite, then there is a (non-canonical) isomorphism with the
$\QQ/\ZZ$-dual $A \cong A^D$, and in particular $|A^D| = |A|$.

\paragraph{Schemes.}
In this text, $X$ always denotes a \textbf{one-dimensional arithmetic scheme},
i.e., a separated scheme of finite type $X \to \Spec \ZZ$ of Krull dimension
$1$.

We remark that the restriction that $X$ is abelian
(Definition~\ref{dfn:abelian-scheme}) is needed only for the
proofs of Theorem~\ref{main-theorem} in \S\ref{sec:direct-proof} and
\S\ref{sec:Weil-etale-proof}. Our calculations in
\S\S\ref{sec:vanishing-order},
\ref{sec:GR-equivariant-cohomology},
\ref{sec:etale-motivic-cohomology},
\ref{sec:regulators},
\ref{sec:Weil-etale-cohomology-of-X}
work for any one-dimensional arithmetic scheme $X$.

\paragraph{Weights.}
In this text, $n$ always stands for a fixed, \emph{strictly negative} integer.

\paragraph{Motivic cohomology.}
We will work with a version of \'{e}tale motivic cohomology defined in terms of
\textbf{Bloch's cycle complexes}. These were introduced by Bloch in
\cite{Bloch-1986} for varieties over fields, and for the version over
$\Spec \ZZ$ see \cite{Geisser-2004-Dedekind,Geisser-2005}.

In short, we let $\Delta^i = \Spec \ZZ [t_0,\ldots,t_i] / (1 - \sum_i t_i)$ be
the algebraic simplex. Denote by $z_n (X, i)$ the group freely generated by
algebraic cycles $Z \subset X \times \Delta^i$ of dimension $n+i$ that intersect
the faces properly. For $n < 0$ we consider the complex of sheaves on $X_\et$
$$\ZZ^c (n) \dfn z_n (\text{\textvisiblespace}, -\bullet) [2n].$$
The corresponding (hyper)cohomology
$$H^i (X_\et, \ZZ^c (n)) \dfn H^i (R\Gamma (X_\et, \ZZ^c (n)))$$
is what we will call in this text \textbf{(\'{e}tale) motivic cohomology}.
For a proper regular arithmetic scheme $X$ of pure dimension $d$ we have
\begin{equation}
  \label{eqn:Zc(n)-vs-Z(d-n)}
  \ZZ^c (n) \cong \ZZ (d-n) [2d],
\end{equation}
where $\ZZ (m)$ is the other motivic complex that usually appears in the
literature; see \cite{Geisser-2004-Dedekind,Geisser-2005} for the definition.
To avoid any confusion, all our calculations will be in terms of $\ZZ^c (n)$.

By \cite[Corollary~7.2]{Geisser-2010}, the groups $H^i (X_\et, \ZZ^c (n))$
satisfy the \textbf{localization property}: if $Z \subset X$ is a closed
subscheme and $U = X\setminus Z$ is its closed complement, then there is a
distinguished triangle
\[ R\Gamma (Z_\et, \ZZ^c (n)) \to
  R\Gamma (X_\et, \ZZ^c (n)) \to
  R\Gamma (U_\et, \ZZ^c (n)) \to 
  R\Gamma (Z_\et, \ZZ^c (n)) [1], \]
giving a long exact sequence
\begin{multline}
  \label{eqn:localization-les}
  \cdots \to H^i (Z_\et, \ZZ^c (n)) \to
  H^i (X_\et, \ZZ^c (n)) \to
  H^i (U_\et, \ZZ^c (n)) \to \\
  H^{i+1} (Z_\et, \ZZ^c (n)) \to \cdots
\end{multline}
This means that $H^i (-, \ZZ^c (n))$ behaves like (motivic) Borel--Moore
homology.

At the level of zeta functions, the localization property corresponds to the
identity
$$\zeta (X,s) = \zeta (Z,s)\,\zeta (U,s).$$
For more results on $\ZZ^c (n)$, we refer the reader to \cite{Geisser-2010}.
 
In general, the groups $H^i (X_\et, \ZZ^c (n))$ are very hard to
compute. However, they are quite well understood for one-dimensional arithmetic
schemes $X$; see \S\ref{sec:etale-motivic-cohomology} below.

\subsection*{Outline of the paper}

In \S\ref{sec:devissage} we prove a d\'{e}vissage lemma that shows how a
property that holds for curves over finite fields and for number rings can be
generalized to any one-dimensional arithmetic scheme. It is an elementary
argument, isolated to avoid repeating the same reasoning in several proofs.

In \S\ref{sec:vanishing-order} we calculate the vanishing order of $\zeta (X,s)$
at $s = n < 0$. Then in \S\ref{sec:GR-equivariant-cohomology} we calculate the
$G_\RR$-equivariant cohomology groups of the finite discrete space of complex
points $X(\CC)$. In \S\ref{sec:etale-motivic-cohomology} we put together various
well-known results to describe the motivic cohomology groups
$H^i (X_\et, \ZZ^c(n))$. In \S\ref{sec:regulators} we define the regulator that
appears in the special value formula.

Our first ``elementary'' proof of the main result is given in
\S\ref{sec:direct-proof}. Then \S\ref{sec:Weil-etale-cohomology-of-X} is devoted
to a calculation of the Weil-\'{e}tale cohomology groups $H^i_\Wc (X, \ZZ(n))$
from \cite{Beshenov-Weil-etale-1} for one-dimensional $X$, which we consider an
interesting result on its own. We use these calculations in
\S\ref{sec:Weil-etale-proof} to formulate explicitly the conjecture
$\mathbf{C} (X,n)$ from \cite{Beshenov-Weil-etale-2}, again for one-dimensional
$X$. This is a second, more conceptual proof of the main result, and it explains
how we arrived at \eqref{eqn:special-value-formula} in the first place.

Finally, we conclude in \S\ref{sec:examples} with a couple of examples showing
how our special value formula works.

\subsection*{Acknowledgments}

I am grateful to Baptiste Morin for various discussions that led to this text.


\section{D\'{e}vissage lemma for one-dimensional schemes}
\label{sec:devissage}

The main idea of this paper is to consider a property that holds for spectra of
number rings $X = \Spec \mathcal{O}_F$ and curves over finite fields $X/\FF_q$,
and then generalize it formally to any one-dimensional arithmetic scheme.
To this end, in this section we isolate a d\'{e}vissage argument which will be
used repeatedly in the rest of the paper.

\begin{lemma}
  \label{lemma:devissage}
  Let $\mathcal{P}$ be a property of arithmetic schemes of Krull dimension
  $\le 1$. Suppose that it satisfies the following compatibilities.
  \begin{enumerate}
  \item[a)] $\mathcal{P} (X)$ holds if and only if $\mathcal{P} (X_\red)$ holds.

  \item[b)] If $X = \coprod_i X_i$ is a finite disjoint union, then
    $\mathcal{P} (X)$ is equivalent to the conjunction of $\mathcal{P} (X_i)$ for
    all $i$.

  \item[c)] If $U \subset X$ is a dense open subscheme, then $\mathcal{P} (X)$
    is equivalent to $\mathcal{P} (U)$.
  \end{enumerate}
  Suppose that
  \begin{enumerate}
  \item[0)] $\mathcal{P} (\Spec \FF_q)$ holds for any finite field $\FF_q$,

  \item[1)] $\mathcal{P} (X)$ holds for any smooth curve $X/\FF_q$,

  \item[2)] $\mathcal{P} (\Spec \mathcal{O}_F)$ holds for any number field
    $F/\QQ$.
  \end{enumerate}
  Then $\mathcal{P} (X)$ holds for any one-dimensional arithmetic scheme $X$.

  \begin{proof}
    First suppose that $\dim X = 0$. Then, thanks to a), we can assume that $X$
    is reduced, and then $X = \coprod_i \Spec \FF_{q,i}$ is a finite disjoint
    union of spectra of finite fields such that $\mathcal{P} (X)$ holds thanks
    to 0) and b).

    Now consider the case of $\dim X = 1$. Again, we can assume that $X$ is
    reduced. We take the normalization $\nu\colon X' \to X$. This is a
    birational morphism: there are dense open subschemes $U \subset X$ and
    $U' \subset X'$ such that
    $\left.\nu\right|_{U'}\colon U' \xrightarrow{\cong} U$ is an
    isomorphism. Thanks to c), we have
    \[ \mathcal{P} (X) \iff
      \mathcal{P} (U) \iff
      \mathcal{P} (U') \iff
      \mathcal{P} (X'). \]
    Therefore, we can assume that $X$ is regular. Now $X = \coprod_i X_i$ is a
    finite disjoint union of normal integral schemes, so thanks to b), we can
    assume that $X$ is integral. There are two cases.

    \begin{itemize}
    \item If $X \to \Spec \ZZ$ lives over a closed point, then it is a smooth
      curve over $\FF_q$, and the claim holds thanks to 1).

    \item If $X \to \Spec \ZZ$ is a dominant morphism, consider an open affine
      neighborhood of the generic point $U \subset X$. Again, $\mathcal{P} (X)$
      is equivalent to $\mathcal{P} (U)$, so it suffices to prove the claim for
      $U$. We have $U = \Spec \mathcal{O}_{F,S}$ for a number field $F/\QQ$ and
      a finite set of places $S$, so everything reduces to
      $\mathcal{P} (\Spec \mathcal{O}_F)$. \qedhere
    \end{itemize}
  \end{proof}
\end{lemma}


\section{Vanishing order of $\zeta (X,s)$ at $s = n < 0$}
\label{sec:vanishing-order}

\begin{definition}[Numbers $r_1$ and $r_2$]
  Given a one-dimensional arithmetic scheme $X$, consider the finite discrete
  space of complex points
  $$X(\CC) \dfn \Hom (\Spec \CC, X).$$
  There is a canonical action of the complex conjugation
  $G_\RR \dfn \Gal (\CC/\RR)$ on $X(\CC)$. The fixed points of this action
  correspond to the real points $X (\RR)$, also known as the
  \textbf{real places}. We set $r_1 = |X (\RR)|$. The non-real places are called
  \textbf{complex places}. They come in conjugate pairs, and we denote their
  number by $2 r_2$.
\end{definition}

\begin{figure}
  \[ \begin{tikzpicture}
      \matrix(m)[matrix of math nodes, row sep=1em, column sep=1em,text height=1ex, text depth=0.2ex]{
        ~ & ~ & ~ & ~ & ~ & \bullet & \bullet & \cdots & \bullet \\
        \bullet & \bullet & \cdots & \bullet \\
        ~ & ~ & ~ & ~ & ~ & \bullet & \bullet & \cdots & \bullet \\};

      \draw[->] (m-2-1) edge[loop above,min distance=10mm] (m-2-1);
      \draw[->] (m-2-2) edge[loop above,min distance=10mm] (m-2-2);
      \draw[->] (m-2-4) edge[loop above,min distance=10mm] (m-2-4);

      \draw[->] (m-1-6) edge[bend left] (m-3-6);
      \draw[->] (m-1-7) edge[bend left] (m-3-7);
      \draw[->] (m-1-9) edge[bend left] (m-3-9);

      \draw[->] (m-3-6) edge[bend left] (m-1-6);
      \draw[->] (m-3-7) edge[bend left] (m-1-7);
      \draw[->] (m-3-9) edge[bend left] (m-1-9);

      \draw [decorate,decoration={brace,amplitude=5pt,mirror}] ($(m-3-1)+(-0.5em,-0.5em)$) -- ($(m-3-4)+(0.5em,-0.5em)$);
      \draw [decorate,decoration={brace,amplitude=5pt,mirror}] ($(m-3-6)+(-0.5em,-0.5em)$) -- ($(m-3-9)+(0.5em,-0.5em)$);

      \draw ($(m-3-1)!.5!(m-3-4)$) node[yshift=-2em,anchor=base] {$r_1$ points};
      \draw ($(m-3-6)!.5!(m-3-9)$) node[yshift=-2em,anchor=base] {$2 r_2$ points};
    \end{tikzpicture} \]

  \caption{$G_\RR \dfn \Gal (\CC/\RR)$ acting on $X (\CC)$}
  \label{fig:complex-conjugation-on-X(C)}
\end{figure}

Equivalently, for a number field $F/\QQ$, denote by $r_1 (F)$ the number of real
embeddings $F \hookrightarrow \RR$ and by $r_2 (F)$ the number of pairs of
complex embeddings $F \hookrightarrow \CC$. Then $r_1 (F) = r_1$ and
$r_2 (F) = r_2$ for $X = \Spec \mathcal{O}_F$. In general, for a one-dimensional
arithmetic scheme $X$, we have
\begin{align*}
  r_1 & = \sum_{\fchar \kappa (\eta) = 0} r_1 (\kappa (\eta)), \\
  r_2 & = \sum_{\fchar \kappa (\eta) = 0} r_2 (\kappa (\eta)),
\end{align*}
where the sums are over generic points $\eta \in X$ with residue field
$\kappa (\eta)$ of characteristic $0$.

\begin{proposition}
  \label{prop:vanishing-order-equals-dn}
  Let $X$ be a one-dimensional arithmetic scheme with $r_1$ real and $2r_2$
  complex places. For $n < 0$, the vanishing order of $\zeta (X,s)$ at $s = n$
  is given by
  \begin{equation}
    \label{eqn:dn}
    d_n = \ord_{s = n} \zeta (X,s) =
    \begin{cases}
      r_1 + r_2, & n\text{ even}, \\
      r_2, & n\text{ odd}.
    \end{cases}
  \end{equation}

  \begin{proof}
    For $X = \Spec \mathcal{O}_F$ the claim is a well-known consequence of the
    functional equation for the Dedekind zeta function
    \cite[\S VII.5]{Neukirch-1999}. It also holds for $X/\FF_q$ since in this
    case $\zeta (X,s)$ has no zeros or poles at $s = n < 0$ according to
    \cite[pp.\,26--27]{Katz-1994}. We now proceed using
    Lemma~\ref{lemma:devissage}.

    We have $\zeta (X,s) = \zeta (X_\red,s)$ and
    $r_{1,2} (X) = r_{1,2} (X_\red)$.
    If $X = \coprod_i X_i$ is a finite disjoint union, then
    \begin{align*}
      \ord_{s = n} \zeta (X,s) & = \sum_i \ord_{s = n} \zeta (X_i,s), \\
      r_{1,2} (X) & = \sum_i r_{1,2} (X_i),
    \end{align*}
    so that the property is compatible with disjoint unions. Finally, if
    $U \subset X$ is a dense open subscheme, then $Z = X\setminus U$ is a
    zero-dimensional scheme, and
    \begin{align*}
      \ord_{s = n} \zeta (X,s) & = \ord_{s = n} \zeta (U,s), \\
      r_{1,2} (X) & = r_{1,2} (U),
    \end{align*}
    so that the property is compatible with taking dense open subschemes.
    We conclude that Lemma~\ref{lemma:devissage} applies.
  \end{proof}
\end{proposition}


\section{$G_\RR$-equivariant cohomology of $X(\CC)$}
\label{sec:GR-equivariant-cohomology}

Viewing $\ZZ (n) \dfn (2\pi i)^n\,\ZZ$ as a constant $G_\RR$-equivariant sheaf
on $X(\CC)$, we consider the $G_\RR$-equivariant cohomology groups (resp. Tate
cohomology)
\begin{align*}
  H^i_c (G_\RR, X(\CC), \ZZ(n)) & \dfn H^i \Bigl(R\Gamma (G_\RR, R\Gamma_c (X(\CC), \ZZ(n)))\Bigr), \\
  \widehat{H}^i_c (G_\RR, X(\CC), \ZZ(n)) & \dfn H^i \Bigl(R\widehat{\Gamma} (G_\RR, R\Gamma_c (X(\CC), \ZZ(n)))\Bigr).
\end{align*}
Of course, $X(\CC)$ is just a finite discrete space, so it is not necessary to
use cohomology with compact support, but we use this notation for consistency
with the general case considered in \cite{Beshenov-Weil-etale-1}. Since
$\dim X(\CC) = 0$, we have
\begin{align*}
  H^i_c (G_\RR, X(\CC), \ZZ(n)) & \cong H^i (G_\RR, H^0_c (X(\CC), \ZZ(n))), \\
  \widehat{H}^i_c (G_\RR, X(\CC), \ZZ(n)) & \cong \widehat{H}^i (G_\RR, H^0_c (X(\CC), \ZZ(n))).
\end{align*}

\begin{proposition}
  Let $X$ be a one-dimensional arithmetic scheme with $r_1$ real places.
  Then the $G_\RR$-equivariant cohomology of $X(\CC)$ is
  \begin{align}
    \label{eqn:Tate-GR-cohomology-of-X(C)}
    \widehat{H}^i_c (G_\RR, X(\CC), \ZZ (n)) & \cong
                                              \begin{cases}
                                                (\ZZ/2\ZZ)^{\oplus r_1}, & i \equiv n ~ (2), \\
                                                0, & i \not\equiv n ~ (2);
                                              \end{cases} \\
    \label{eqn:usual-GR-cohomology-of-X(C)}
    H^i_c (G_\RR, X(\CC), \ZZ (n)) & \cong
                                    \begin{cases}
                                      0, & i < 0, \\
                                      \ZZ^{\oplus d_n}, & i = 0, \\
                                      \widehat{H}^i_c (G_\RR, X(\CC), \ZZ (n)), & i \ge 1.
                                    \end{cases}
  \end{align}
  Here $d_n$ is the vanishing order given by \eqref{eqn:dn}.

  \begin{proof}
    We have
    \[ H^0_c (X(\CC), \ZZ(n)) \cong
      \ZZ (n)^{\oplus r_1} \oplus (\ZZ (n) \oplus \ZZ (n))^{\oplus r_2}, \]
    and the $G_\RR$-action on the two summands is given by
    $x \mapsto \overline{x}$ and $(x,y) \mapsto (\overline{y}, \overline{x})$,
    respectively. (See Figure~\eqref{fig:complex-conjugation-on-X(C)}.)

    We recall that the Tate cohomology of a finite cyclic group is $2$-periodic:
    \[ \widehat{H}^i (G,A) \cong
      \begin{cases}
        \widehat{H}^0 (G,A), & i\text{ even}, \\
        \widehat{H}_0 (G,A), & i\text{ odd},
      \end{cases} \]
    and the groups $\widehat{H}^0 (G,A)$ and $\widehat{H}_0 (G,A)$ are given by
    the exact sequence
    \[ 0 \to \widehat{H}_0 (G,A) \to
      A_G \xrightarrow{N} A^G \to
      \widehat{H}^0 (G,A) \to 0 \]
    where $N$ is the norm map induced by the action of $\sum_{g\in G} g$.

    \vspace{1em}

    Therefore, we can consider two cases.
    \begin{enumerate}
    \item[1)] \textbf{$G_\RR$-module $A = \ZZ (n)$ with the action via
        $x \mapsto \overline{x}$.}

      In this case, we see that
      \[ A^{G_\RR} \cong
        \begin{cases}
          \ZZ, & n\text{ even}, \\
          0, & n\text{ odd}.
        \end{cases} \]
      Similarly, it is straightforward to calculate the coinvariants
      $A_{G_\RR}$, and
      \[ \widehat{H}^0 (G_\RR, A) \cong
        \begin{cases}
          \ZZ/2\ZZ, & n\text{ even},\\
          0, & n\text{ odd},
        \end{cases} \quad
        \widehat{H}_0 (G_\RR, A) \cong
        \begin{cases}
          0, & n\text{ even},\\
          \ZZ/2\ZZ, & n\text{ odd}.
        \end{cases} \]

    \item[2)] \textbf{$G_\RR$-module $A = \ZZ (n) \oplus \ZZ (n)$ with the
        action via $(x,y) \mapsto (\overline{y}, \overline{x})$}.

      In this case $A^{G_\RR} \cong \ZZ$ and
      $\widehat{H}^0 (G_\RR,A) = \widehat{H}_0 (G_\RR,A) = 0$.
    \end{enumerate}

    Combining these two calculations, we obtain Tate cohomology groups
    \eqref{eqn:Tate-GR-cohomology-of-X(C)}. For the usual cohomology
    \eqref{eqn:usual-GR-cohomology-of-X(C)}, we have
    \begin{align*}
      H^0_c (G_\RR, X(\CC), \ZZ (n)) & \cong H^0_c (X(\CC), \ZZ (n))^{G_\RR}, \\
      H^i_c (G_\RR, X(\CC), \ZZ (n)) & \cong \widehat{H}^i_c (G_\RR, X(\CC), \ZZ (n)) \quad \text{for }i \ge 1. \qedhere
    \end{align*}
  \end{proof}
\end{proposition}


\section{\'{E}tale motivic cohomology of one-dimensional schemes}
\label{sec:etale-motivic-cohomology}

In this section we review the structure of the \'{e}tale motivic cohomology
$H^i (X_\et, \ZZ^c(n))$ for one-dimensional $X$ and $n < 0$. What follows is
fairly well-known, so we claim no originality here, but we compile the
references and state the result for a general one-dimensional arithmetic scheme.

\begin{proposition}
  \label{prop:structure-of-motivic-cohomology}
  If $X$ is a one-dimensional arithmetic scheme and $n < 0$, then
  \begin{equation}
    \label{eqn:structure-of-motivic-cohomology}
    H^i (X_\et, \ZZ^c (n)) \cong
    \begin{cases}
      0, & i < -1, \\
      \text{finitely generated of rk } d_n, & i = -1, \\
      \text{finite}, & i = 0,1, \\
      (\ZZ/2\ZZ)^{\oplus r_1}, & i \ge 2, ~ i\not\equiv n ~ (2), \\
      0, & i \ge 2, ~ i\equiv n ~ (2).
    \end{cases}
  \end{equation}
  Here $d_n$ is given by \eqref{eqn:dn} and $r_1 = |X (\RR)|$ is the number of
  real places of $X$. Further, if $X = \Spec \mathcal{O}_F$ for a number field
  $F/\QQ$, then
  \begin{equation}
    \label{eqn:2-torsion-in-H1-Zc-for-Spec-OF}
    H^1 (X_\et, \ZZ^c (n)) \cong
    \begin{cases}
      (\ZZ/2\ZZ)^{\oplus r_1}, & n\text{ even}, \\
      0, & n \text{ odd}.
    \end{cases}
  \end{equation}
\end{proposition}

An important ingredient of our proof is the arithmetic duality
\cite[Theorem~I]{Beshenov-Weil-etale-1}, which states that if
$H^i (X_\et, \ZZ^c (n))$ are finitely generated groups for all $i \in \ZZ$, then
\begin{equation}
  \label{eqn:arithmetic-duality}
  \widehat{H}^i_c (X_\et, \ZZ (n)) \cong
  H^{2-i} (X_\et, \ZZ^c (n))^D,
\end{equation}
where
\begin{equation}
  \label{eqn:definition-of-Z(n)}
  \ZZ (n) \dfn \QQ/\ZZ (n) [-1] \dfn
  \bigoplus_p \varinjlim_r j_{p!} \mu_{p^r}^{\otimes n} [-1].
\end{equation}
Here $\widehat{H}^i_c (X_\et, \ZZ (n))$ is the modified cohomology with compact
support, for which we refer to \cite[\S 2]{Geisser-Schmidt-2018} and
\cite[Appendix~B]{Beshenov-Weil-etale-1}. In particular,
\[
  \widehat{H}^i_c (X_\et, \ZZ (n)) = H^i_c (X_\et, \ZZ (n))
  \quad
  \text{if }X (\RR) = \emptyset.
\]
We recall that $(-)^D$ denotes the group $\Hom (-, \QQ/\ZZ)$. We note that
\eqref{eqn:arithmetic-duality} is a powerful result, deduced in
\cite{Beshenov-Weil-etale-1} from the work of Geisser \cite{Geisser-2010}.

\begin{proof}[Proof of Proposition~\ref{prop:structure-of-motivic-cohomology}]
  We use Lemma~\ref{lemma:devissage}. We will say that $\mathcal{P} (X)$ holds
  if the motivic cohomology of $X$ has the structure
  \eqref{eqn:structure-of-motivic-cohomology}.

  \vspace{1em}

  \textbf{Let us first consider the case of a finite field $X = \Spec \FF_q$}.
  We have
  \begin{equation}
    \label{eqn:motivic-cohomology-finite-fields}
    H^i (\Spec \FF_{q,\et}, \ZZ^c (n)) \cong
    \begin{cases}
      \ZZ/(q^{-n} - 1), & i = 1, \\
      0, & i \ne 1.
    \end{cases}
  \end{equation}
  ---see, for example, \cite[Example~4.2]{Geisser-2017}. This is related to
  Quillen's calculation of the $K$-theory of finite fields \cite{Quillen-1972}.

  In general, if $X$ is a zero-dimensional arithmetic scheme, then the motivic
  cohomology of $X$ and $X_\red$ coincide, so we can assume that $X$ is
  reduced. Then $X$ is a finite disjoint union of $X_i = \Spec \FF_{q_i}$, and
  \begin{equation}
    H^i (X, \ZZ^c (n)) = \begin{cases}
      \text{finite}, & i = 1, \\
      0, & i \ne 1.
    \end{cases}
  \end{equation}
  In particular, $\mathcal{P} (X)$ holds if $\dim X = 0$.

  \vspace{1em}

  \textbf{Now we check the compatibility properties for $\mathcal{P}$}.
  If $X = \coprod_i X_i$ is a finite disjoint union, then
  $H^i (X_\et, \ZZ^c (n)) \cong \bigoplus_i H^i (X_{i,\et}, \ZZ^c (n))$,
  hence the property $\mathcal{P}$ is compatible with disjoint unions.

  Similarly, if $U \subset X$ is a dense open subscheme, and $Z = X\setminus U$
  its closed complement, then $\dim Z = 0$. We consider the long exact sequence
  \eqref{eqn:localization-les}. Since the cohomology of $Z$ is concentrated in
  $i = 1$, we have $H^i (X_\et, \ZZ^c (n)) \cong H^i (U_\et, \ZZ^c (n))$ for
  $i \ne 0,1$, and what is left is an exact sequence
  \begin{multline*}
    0 \to H^0 (X_\et, \ZZ^c (n)) \to
    H^0 (U_\et, \ZZ^c (n)) \to \\
    H^1 (Z_\et, \ZZ^c (n)) \to
    H^1 (X_\et, \ZZ^c (n)) \to
    H^1 (U_\et, \ZZ^c (n)) \to 0
  \end{multline*}
  Moreover, $d_n (X) = d_n (U)$. These considerations show that
  $\mathcal{P} (X)$ and $\mathcal{P} (U)$ are equivalent, and therefore
  Lemma~\ref{lemma:devissage} works, and it remains to establish
  $\mathcal{P} (X)$ for a curve $X/\FF_q$ or $X = \Spec \mathcal{O}_F$.

  \vspace{1em}

  \textbf{Suppose that $X/\FF_q$ is a smooth curve}. The groups
  $H^i (X_\et, \ZZ^c (n))$ are finitely generated by
  \cite[Proposition~4.3]{Geisser-2017}, so that the duality
  \eqref{eqn:arithmetic-duality} holds. The $\QQ/\ZZ$-dual groups
  \[ H^i_c (X_\et, \ZZ(n)) =
    \bigoplus_\ell H^{i-1}_c (X_\et, \QQ_\ell/\ZZ_\ell (n)) \]
  are finite by \cite[Theorem~3]{Kahn-2003}, and concentrated in $i = 1,2,3$ for
  dimension reasons. It follows that $H^i (X_\et, \ZZ^c (n))$ in this case are
  finite groups concentrated in $i = -1,0,1$, and the property $\mathcal{P} (X)$
  holds.

  \vspace{1em}

  \textbf{It remains to consider the case of $X = \Spec \mathcal{O}_F$}.
  In this case, the finite generation of $H^i (X_\et, \ZZ^c (n))$ is also known;
  see, for example, \cite[Proposition~4.14]{Geisser-2017}. Therefore, the
  duality \eqref{eqn:arithmetic-duality} holds. We have
  $\widehat{H}^i_c (\Spec \mathcal{O}_F [1/p], \mu_{p^r}^{\otimes n}) = 0$ for
  $i \ge 3$ by Artin--Verdier duality \cite[Chapter~II,
  Corollary~3.3]{Milne-ADT}, or by \cite[p.\,268]{Soule-1979}. Therefore, it
  follows that $\widehat{H}^i_c (X_\et, \ZZ (n)) = 0$ for $i \ge 4$, and hence
  by duality \eqref{eqn:arithmetic-duality}, $H^i (X_\et, \ZZ^c (n)) = 0$ for
  $i \le -2$.

  Now we identify the finite $2$-torsion in $H^i (X_\et, \ZZ^c (n))$ for
  $i \ge 2$. By \cite[Lemma~6.14]{Flach-Morin-2018}, there is an exact
  sequence
  \begin{multline}
    \label{eqn:les-Hc-vs-H-hat-c}
    \cdots \to H^{i-1}_c (X_\et, \ZZ (n)) \to
    \widehat{H}^{i-1} (G_\RR, X(\CC), \ZZ (n)) \to \\
    \widehat{H}^i_c (X_\et, \ZZ (n)) \to
    H^i_c (X_\et, \ZZ (n)) \to \cdots
  \end{multline}
  For $i \le 0$ we have $H^i_c (X_\et, \ZZ (n)) = 0$, and therefore
  \[ \widehat{H}^i_c (X_\et, \ZZ (n)) \cong
    \widehat{H}^{i-1}_c (G_\RR, X(\CC), \ZZ (n)) \cong
    \begin{cases}
      (\ZZ/2\ZZ)^{\oplus r_1}, & i\not\equiv n ~ (2), \\
      0, & i\equiv n ~ (2).
    \end{cases} \]
  By duality, for $i \ge 2$ we have
  \[ H^i (X_\et, \ZZ^c (n)) \cong
    \begin{cases}
      (\ZZ/2\ZZ)^{\oplus r_1}, & i\not\equiv n ~ (2), \\
      0, & i\equiv n ~ (2).
    \end{cases} \]

  Now we determine the ranks of $H^i (X_\et, \ZZ^c (n))$ for $i = -1,0,1$.
  By \cite[Proposition~2.1]{Kolster-Sands-2008} the Chern character for
  $i = -1,0$
  $$K_{-2n - i} (X) \to H^i (X_\et, \ZZ^c(n))$$
  has a finite $2$-torsion kernel and cokernel.
  (Originally, the target group is defined over $X_\Zar$, and we identify it
  with the cohomology on $X_\et$ using the Beilinson--Lichtenbaum conjecture,
  which is now a theorem \cite[Theorem~1.2]{Geisser-2004-Dedekind}. We further
  use the isomorphism \eqref{eqn:Zc(n)-vs-Z(d-n)} to identify our motivic
  cohomology with the one used in \cite{Kolster-Sands-2008}.)

  For $i = -1,0$ we have therefore
  $$\rk_\ZZ H^i (X_\et, \ZZ^c(n)) = \rk_\ZZ K_{-2n - i} (X).$$
  Together with Borel's calculation of the ranks of $K_m (\mathcal{O}_F)$ in \cite{Borel-1974}, this implies that
  $H^0 (X_\et, \ZZ^c (n))$ is a finite group, while
  \[ \rk_\ZZ H^{-1} (X_\et, \ZZ^c(n)) = d_n =
    \begin{cases}
      r_1 + r_2, & n \text{ even}, \\
      r_2, & n \text{ odd}.
    \end{cases} \]
  Finally, by \cite[p.\,179]{Kolster-Sands-2008} and
  \eqref{eqn:Zc(n)-vs-Z(d-n)}, we have
  \[ H^1 (X_\et, \ZZ^c(n)) \cong
    \begin{cases}
      (\ZZ/2\ZZ)^{\oplus r_1}, & n\text{ even}, \\
      0, & n\text{ odd}.
    \end{cases} \]
  This concludes the proof.
\end{proof}


\section{Regulator for one-dimensional $X$}
\label{sec:regulators}

Now we explain what is meant by the regulator in our situation.

\begin{definition}
  We let the \textbf{regulator morphism} be the composition
  \begin{multline*}
    \varrho_{X,n}\colon
    H^{-1} (X_\et, \ZZ^c (n)) \xrightarrow{x \mapsto x\otimes 1}
    H^{-1} (X_\et, \ZZ^c (n)) \otimes \RR \\
    \xrightarrow{Reg_{X,n}} H^0_\BM (G_\RR, X(\CC), \RR(n)),
  \end{multline*}
  where the map $Reg_{X,n}$ is defined in \cite[\S 2]{Beshenov-Weil-etale-2}.
\end{definition}

The target is the Borel--Moore cohomology defined by
\[ H^0_\BM (G_\RR, X(\CC), \RR(n)) \dfn
  \Hom (H^0_c (G_\RR, X(\CC), \RR(n)), \RR). \]
In general, the regulator takes values in Deligne--Beilinson cohomology, but the
target simplifies in the case of $n < 0$, as explained in
\cite[\S 2]{Beshenov-Weil-etale-2}.

\begin{remark}
  The only relevant group for the regulator is $H^{-1} (X_\et, \ZZ^c (n))$,
  since the cohomology in other degrees is finite by
  Proposition~\ref{prop:structure-of-motivic-cohomology}.

  The general definition in \cite[\S 2]{Beshenov-Weil-etale-2} is based on the
  construction of Kerr, Lewis and M\"{u}ller-Stach
  \cite{Kerr-Lewis-Muller-Stach-2006} which works at the level of
  complexes. This is not very important in the one-dimensional case, where the
  interesting cohomology is concentrated in $i = -1$. The reader can use any
  other equivalent construction of the regulator for motivic cohomology.
\end{remark}

\begin{remark}
  If $X = \Spec \mathcal{O}_F$, then $\varrho_{X,n}$ can be identified with the
  Beilinson regulator map that appears in the special value conjecture of Flach
  and Morin in \cite[\S 5.8.3]{Flach-Morin-2018}.
\end{remark}

\begin{lemma}
  \label{lemma:regulator-isomorphism}
  For any one-dimensional arithmetic scheme $X$ and $n < 0$, the $\RR$-dual to
  the regulator
  \[ Reg_{X,n}^\vee\colon H^0_c (G_\RR, X(\CC), \RR(n)) \to
    \Hom (H^{-1} (X_\et, \ZZ^c (n)), \RR) \]
  is an isomorphism.

  \begin{proof}
    If $X/\FF_q$, then the claim is trivial. For $X = \Spec \mathcal{O}_F$, this
    is a well-known property of the Beilinson regulator. To apply
    Lemma~\ref{lemma:devissage}, we need to check compatibility with disjoint
    unions and passing to a dense open subscheme $U \subset X$. For disjoint
    unions, this is clear. For a dense open subscheme $U \subset X$, the closed
    complement $Z = X\setminus U$ has dimension $0$, and the localization exact
    sequence \eqref{eqn:localization-les} with the long exact sequence for
    cohomology with compact support yields integral isomorphisms
    \begin{align*}
      H^{-1} (X_\et, \ZZ^c (n)) & \xrightarrow{\cong} H^{-1} (U_\et, \ZZ^c (n)), \\
      H^0_c (G_\RR, U(\CC), \ZZ(n)) & \xrightarrow{\cong} H^0_c (G_\RR, Z(\CC), \ZZ(n)).
    \end{align*}

    We now have a commutative diagram
    \[ \begin{tikzcd}
        H^0_c (G_\RR, U(\CC), \RR(n)) \ar{r}{Reg_{U,n}^\vee} \ar{d}{\cong} &
        \Hom (H^{-1} (U_\et, \ZZ^c (n)), \RR) \ar{d}{\cong} \\
        H^0_c (G_\RR, X(\CC), \RR(n)) \ar{r}{Reg_{X,n}^\vee} &
        \Hom (H^{-1} (X_\et, \ZZ^c (n)), \RR)
      \end{tikzcd} \]
    The upper arrow is an isomorphism if and only if the lower arrow is.
  \end{proof}
\end{lemma}

\begin{definition}
  \label{dfn:regulator}
  For a one-dimensional arithmetic scheme $X$, we define the \textbf{regulator}
  to be
  \[ R_{X,n} \dfn \vol \Bigl(\coker \bigl(
    H^{-1} (X_\et, \ZZ^c (n)) \xrightarrow{\varrho_{X,n}}
    H^0_\BM (G_\RR, X(\CC), \RR(n))
    \bigr)\Bigr), \]
  where the volume is taken with respect to the canonical integral structure.
\end{definition}

If $X(\CC) = \emptyset$, or $n$ is odd and $r_2 = 0$, then
$H^0_\BM (G_\RR, X(\CC), \RR (n)) = 0$, and we set $R_{X,n} = 1$.

\begin{lemma}
  \label{lemma:regulator-dense-open-subset}
  Let $X$ be a one-dimensional arithmetic scheme and $n < 0$. For any dense open
  subscheme $U \subset X$, we have $R_{X,n} = R_{U,n}$.

  \begin{proof}
    Follows from the proof of Lemma~\ref{lemma:regulator-isomorphism}.
  \end{proof}
\end{lemma}

\begin{proposition}
  \label{prop:trivialization-of-free-part}
  Given a one-dimensional arithmetic scheme $X$ and $n < 0$, consider the
  two-term acyclic complex of real vector spaces
  \[ C^\bullet\colon
    0 \to
    \underbrace{H^0_c (G_\RR, X(\CC), \RR(n))}_{\deg 0}
    \xrightarrow{Reg_{X,n}^\vee}
    \underbrace{\Hom (H^{-1} (X_\et, \ZZ^c (n)), \RR)}_{\deg 1}
    \to 0 \]
  Then taking the determinant $\det_\RR (C^\bullet)$ in the sense of
  Knudsen and Mumford \cite{Knudsen-Mumford-1976},
  the image of the canonical map
  \begin{multline*}
    \det_\ZZ H^0_c (G_\RR, X(\CC), \ZZ(n)) \otimes_\ZZ
    \det_\ZZ \Hom (H^{-1} (X_\et, \ZZ^c (n)), \ZZ)^{-1} \to \\
    \det_\RR H^0_c (G_\RR, X(\CC), \RR (n)) \otimes_\RR
    \det_\RR \Hom (H^{-1} (X_\et, \ZZ^c (n)), \RR)^{-1}
    \xrightarrow{\cong} \RR
  \end{multline*}
  corresponds to $R_{X,n}\,\ZZ \subset \RR$.

  \begin{proof}
    In general, if $F$ and $F'$ are free groups of finite rank $d$, and
    $$C^\bullet\colon 0 \to F\otimes_\ZZ \RR \xrightarrow{\phi} F'\otimes_\ZZ \RR \to 0$$
    is a two-term acyclic complex of real vector spaces, then the image of
    \[ \ZZ \cong \det_\ZZ F \otimes_\ZZ (\det_\ZZ F')^{-1} \to
      \det_\RR (F\otimes_\ZZ \RR) \otimes_\RR \det_\RR (F' \otimes_\ZZ \RR)^{-1}
      = \det_\RR (C^\bullet) \xrightarrow{\cong} \RR \]
    corresponds to $D\ZZ \subset \RR$, where $D$ is the determinant of $\phi$ with
    respect to the bases induced by $\ZZ$-bases of $F$ and $F'$.
    This follows from the explicit description of the isomorphism
    $\det_\RR (C^\bullet) \xrightarrow{\cong} \RR$ from
    \cite[p.\,33]{Knudsen-Mumford-1976}: it is
    \[ \det_\RR (F \otimes_\ZZ \RR) \otimes_\RR
      \det_\RR (F'\otimes_\ZZ \RR)^{-1} \xrightarrow{\det_\RR (\phi)}
      \det_\RR (F' \otimes_\ZZ \RR) \otimes_\RR
      \det_\RR (F' \otimes_\ZZ \RR)^{-1} \xrightarrow{\cong} \RR \]
    where the last arrow is the canonical pairing.

    Therefore, in our situation, the image of
    \[ \det_\ZZ H^0_c (G_\RR, X(\CC), \ZZ(n)) \otimes_\ZZ
      \det_\ZZ \Hom (H^{-1} (X_\et, \ZZ^c (n)), \ZZ)^{-1} \]
    is $D \ZZ \subset \RR$, where $D$ is the determinant of $Reg_{X,n}^\vee$
    considered with respect to the bases induced by $\ZZ$-bases of
    $H^0_c (G_\RR, X(\CC), \ZZ(n))$ and
    $\Hom (H^{-1} (X_\et, \ZZ^c (n)), \ZZ)$. Dually, $D = R_{X,n}$.
  \end{proof}
\end{proposition}


\section{Direct proof of the special value formula}
\label{sec:direct-proof}

In this section we explain how to prove our special value formula directly by
combining the known special value formulas for $X = \Spec \mathcal{O}_F$ and
curves over finite fields $X/\FF_q$ via localization.

\begin{lemma}
  \label{lemma:elementary-proof-1}
  Let $n < 0$.

  \begin{enumerate}
  \item[0)] If $X$ is a zero-dimensional arithmetic scheme, then
    $$\zeta (X,n) = \pm \frac{1}{|H^1 (X_\et, \ZZ^c (n))|}.$$

  \item[1)] If $X/\FF_q$ is a curve over a finite field, then
    \[ \zeta (X,n) =
      \pm\frac{|H^0 (X_\et, \ZZ^c(n))|}{|H^{-1} (X_\et, \ZZ^c(n))|\cdot |H^1 (X_\et, \ZZ^c(n))|}. \]

  \item[2)] If $X = \Spec \mathcal{O}_F$ for an abelian number field $F/\QQ$,
    then
    \[ \zeta (X,n) = \pm\frac{|H^0 (X_\et, \ZZ^c(n))|}{|H^{-1} (X_\et, \ZZ^c(n))|}\,R_{X,n}. \]
  \end{enumerate}

  In particular, formula \eqref{eqn:special-value-formula} holds in these
  cases.

  \begin{proof}
    In part 0), motivic cohomology and the zeta function do not distinguish
    between $X$ and $X_\red$, so we can assume that $X$ is a finite   
    disjoint union of $\Spec \FF_{q_i}$. Thanks to
    \eqref{eqn:motivic-cohomology-finite-fields},
    \[ \zeta (X,n) = \prod_i \frac{1}{1 - q_i^{-n}} =
      \pm \prod_i \frac{1}{|H^1 (X_{i,\et}, \ZZ^c (n))|} =
      \pm \frac{1}{|H^1 (X_\et, \ZZ^c (n))|}. \]
    Note that this is formula \eqref{eqn:special-value-formula}
    since $\delta = 0$ in this case and
    $H^{-1} (X_\et, \ZZ^c (n)) = H^0 (X_\et, \ZZ^c (n)) = 0$ by
    \eqref{eqn:motivic-cohomology-finite-fields}.

    \vspace{1em}

    For part 1), we refer the reader to \cite[\S 5]{Beshenov-Weil-etale-2}.
    Part 2) follows from \cite[Proposition~5.35]{Flach-Morin-2018}.
    The formula is equivalent to \eqref{eqn:special-value-formula}, since
    $2^\delta = |H^1 (X_\et, \ZZ (n))|$ by
    \eqref{eqn:2-torsion-in-H1-Zc-for-Spec-OF}.
  \end{proof}
\end{lemma}

\begin{remark}
  The special value at $s = 0$ is not necessarily a rational number:
  \[ \zeta^* (\Spec \FF_q, 0) =
    \lim_{s \to 0} \frac{s}{1 - q^{-s}} =
    \frac{1}{\log q}. \]
  Moreover,
  \[ H^i (\Spec \FF_{q,\et}, \ZZ^c (0)) =
    \begin{cases}
      \ZZ, & i = 1, \\
      \QQ/\ZZ, & i = 2, \\
      0, & i \ne 1,2.
    \end{cases} \]
  This toy example already shows that it is important that we focus on the case
  of $n < 0$.
\end{remark}

\begin{lemma}
  \label{lemma:elementary-proof-2}
  Let $X$ be a one-dimensional arithmetic scheme and let $U \subset X$ be a
  dense open subscheme. Then the special value formula
  \eqref{eqn:special-value-formula} for $X$ is equivalent to the corresponding
  formula for $U$.

  \begin{proof}
    Let $Z = X\setminus U$ be the zero-dimensional complement. We have
    $$\zeta (X,n) = \zeta (U,n)\,\zeta (Z,n),$$
    where
    \begin{align}
      \label{eqn:special-value-for-X} \zeta (X,n) & \stackrel{?}{=} \pm 2^\delta\,\frac{|H^0 (X_\et, \ZZ^c (n)|}{|H^{-1} (X_\et, \ZZ^c (n))_\tors| \cdot |H^1 (X_\et, \ZZ^c (n))|}\,R, \\
      \label{eqn:special-value-for-U} \zeta (U,n) & \stackrel{?}{=} \pm 2^\delta\,\frac{|H^0 (U_\et, \ZZ^c (n)|}{|H^{-1} (U_\et, \ZZ^c (n))_\tors| \cdot |H^1 (U_\et, \ZZ^c (n))|}\,R, \\
      \notag \zeta (Z,n) & = \pm\frac{1}{|H^1 (Z_\et, \ZZ^c (n))|}.
    \end{align}
    Here $\delta = \delta_{X,n} = \delta_{U,n}$, and $R = R_{X,n} = R_{U,n}$
    (see Lemma~\ref{lemma:regulator-dense-open-subset}). We note that
    $|H^{-1} (X_\et, \ZZ^c (n))_\tors| = |H^{-1} (U_\et, \ZZ^c (n))_\tors|$.
    On the other hand, the exact sequence of finite groups
    \begin{multline}
      0 \to H^0 (X_\et, \ZZ^c (n)) \to
      H^0 (U_\et, \ZZ^c (n)) \to \\
      H^1 (Z_\et, \ZZ^c (n)) \to
      H^1 (X_\et, \ZZ^c (n)) \to
      H^1 (U_\et, \ZZ^c (n)) \to 0
    \end{multline}
    gives
    \[ \frac{|H^0 (X_\et, \ZZ^c (n))|}{|H^1 (X_\et, \ZZ^c (n))|} =
      \frac{|H^0 (U_\et, \ZZ^c (n))|}{|H^1 (U_\et, \ZZ^c (n))|}\cdot
      \frac{1}{|H^1 (Z_\et, \ZZ^c (n))|}. \]
    From this we see that \eqref{eqn:special-value-for-X} and
    \eqref{eqn:special-value-for-U} are equivalent.
  \end{proof}
\end{lemma}

The above Lemmas~\ref{lemma:elementary-proof-1} and
\ref{lemma:elementary-proof-2} together with Lemma~\ref{lemma:devissage} now
prove Theorem~\ref{main-theorem} from the introduction.

\begin{remark}
  \label{rmk:TNC}
  Our proof of Lemma~\ref{lemma:elementary-proof-1} uses
  \cite[Proposition~5.35]{Flach-Morin-2018}, which in turn reduces to the
  Tamagawa number conjecture for abelian $F/\QQ$. The non-abelian version of
  Theorem~\ref{main-theorem} is therefore equivalent to the corresponding
  conjecture for non-abelian $F/\QQ$.
\end{remark}

\begin{remark}
  Note that $\zeta (\Spec \FF_q, n) = \frac{1}{1 - q^{-n}} < 0$. Thus, if we
  remove $m$ closed points from $X$, the sign of $\zeta^* (X,n)$ changes by
  $(-1)^m$. It is not hard to figure out the sign in any concrete example;
  however, it is not so clear in what terms to write the general expression for
  the sign.
\end{remark}


\section{Weil-\'{e}tale cohomology of one-dimensional arithmetic schemes}
\label{sec:Weil-etale-cohomology-of-X}

In this section we calculate Weil-\'{e}tale cohomology groups
$H^i_\Wc (X, \ZZ(n))$ for ${n < 0}$, as defined in
\cite{Beshenov-Weil-etale-1}. Let us briefly recall the construction.
In general, let $X$ be an arithmetic scheme with finitely generated motivic
cohomology $H^i (X_\et, \ZZ^c(n))$. The construction is carried out in two
steps.

\begin{itemize}
\item \textbf{Step 1}. Consider the morphism in the derived category
  $\mathbf{D} (\ZZ)$
  \[ \alpha_{X,n}\colon \RHom (R\Gamma (X_\et, \ZZ^c(n)), \QQ [-2]) \to
    R\Gamma_c (X_\et, \ZZ(n)) \]
  determined at the level of cohomology, using the arithmetic duality
  \eqref{eqn:arithmetic-duality}, by
  \begin{multline}
    \label{eqn:alpha-X-n-cohomology}
    H^i (\alpha_{X,n})\colon \Hom (H^{2-i} (X_\et, \ZZ^c(n)), \QQ)
    \xrightarrow{\QQ \twoheadrightarrow \QQ/\ZZ}
    H^{2-i} (X_\et, \ZZ^c(n))^D \\
    \xleftarrow{\cong} \widehat{H}^i_c (X_\et, \ZZ(n)) \to
    H^i_c (X_\et, \ZZ(n)).
  \end{multline}
  The complex $R\Gamma_\fg (X, \ZZ(n))$ is defined as a cone of
  $\alpha_{X,n}$:
  \begin{multline*}
    \RHom (R\Gamma (X_\et, \ZZ^c(n)), \QQ [-2]) \xrightarrow{\alpha_{X,n}}
    R\Gamma_c (X_\et, \ZZ(n)) \\
    \to R\Gamma_\fg (X_\et, \ZZ(n)) \to \RHom (R\Gamma (X_\et, \ZZ^c(n)), \QQ [-1])
  \end{multline*}
  The groups
  $$H^i_\fg (X, \ZZ(n)) \dfn H^i (R\Gamma_\fg (X, \ZZ(n)))$$
  are finitely generated for all $i \in \ZZ$, vanish for $i \ll 0$, and finite
  $2$-torsion for $i \gg 0$. For the details we refer to
  \cite[\S 5]{Beshenov-Weil-etale-1}.

\item \textbf{Step 2}. We consider a canonical morphism $i^*_\infty$ in the
  derived category $\mathbf{D} (\ZZ)$ which is torsion and yields a commutative
  diagram
  \[ \begin{tikzcd}
      R\Gamma_c (X_\et, \ZZ(n)) \ar{r}\ar{d}[swap]{u^*_\infty} & R\Gamma_\fg (X, \ZZ(n)) \ar{dl}{i^*_\infty} \\
      R\Gamma_c (G_\RR, X(\CC), \ZZ(n))
    \end{tikzcd} \]
  ---see \cite[\S\S 6,7]{Beshenov-Weil-etale-1} for more details.
  Weil-\'{e}tale cohomology with compact support is defined as a mapping fiber of
  $i^*_\infty$:
  \[ R\Gamma_\Wc (X, \ZZ(n)) \to
    R\Gamma_\fg (X, \ZZ(n)) \xrightarrow{i^*_\infty}
    R\Gamma_c (G_\RR, X(\CC), \ZZ(n)) \to [1] \]
  The resulting groups
  $$H^i_\Wc (X, \ZZ(n)) \dfn H^i (R\Gamma_\Wc (X, \ZZ(n)))$$
  are finitely generated and vanish for $i \notin [0, 2\dim X + 1]$.
  We refer to \cite[\S 7]{Beshenov-Weil-etale-1} for the general properties.
\end{itemize}

Here we calculate $H^i_\Wc (X, \ZZ(n))$ for one-dimensional $X$.

\begin{proposition}
  \label{prop:calculation-of-H-Wc}
  Let $X$ be a one-dimensional arithmetic scheme and $n < 0$.

  \begin{enumerate}
  \item[0)] $H^i_\Wc (X, \ZZ(n)) = 0$ for $i \ne 1,2,3$.

  \item[1)] There is a short exact sequence
    \begin{equation}
      \label{eqn:H1Wc-ses}
      0 \to \underbrace{H^0_c (G_\RR, X(\CC), \ZZ(n))}_{\cong \ZZ^{\oplus d_n}} \to
      H^1_\Wc (X, \ZZ(n)) \to T_1 \to 0
    \end{equation}
    in which $T_1$ sits in a short exact sequence of finite groups
    \[ 0 \to \widehat{H}^0_c (G_\RR, X(\CC), \ZZ(n)) \to
      H^1 (X_\et, \ZZ^c(n))^D \to
      T_1 \to 0 \]
    In particular, $H^1_\Wc (X, \ZZ(n))$ is finitely generated of rank $d_n$,
    and
    $$|T_1| = \frac{1}{2^\delta}\cdot |H^1 (X_\et, \ZZ^c (n))|,$$
    where $\delta$ is defined by \eqref{eqn:delta}.

  \item[2)] There is an isomorphism of finitely generated groups
    \[ H^2_\Wc (X, \ZZ(n)) \cong
      \underbrace{H^{-1} (X_\et, \ZZ^c(n))^*}_{\cong \ZZ^{\oplus d_n}}
      \oplus
      \underbrace{H^0 (X_\et, \ZZ^c(n))^D}_{\text{finite}}. \]

  \item[3)] There is an isomorphism of finite groups
    $$H^3_\Wc (X, \ZZ(n)) \cong (H^{-1} (X_\et, \ZZ^c(n))_\tors)^D.$$
  \end{enumerate}
\end{proposition}

We recall that $A^D \dfn \Hom (A, \QQ/\ZZ)$ and $A^* \dfn \Hom (A, \ZZ)$.

\begin{proof}
  From the definition of $R\Gamma_\fg (X, \ZZ(n))$ we have a long exact
  sequence
  \begin{multline}
    \label{eqn:fg-les}
    \cdots \to \Hom (H^{2-i} (X_\et, \ZZ^c (n)), \QQ) \xrightarrow{H^i (\alpha_{X,n})} H^i_c (X_\et, \ZZ(n)) \\
    \to H^i_\fg (X, \ZZ(n)) \to \Hom (H^{1-i} (X_\et, \ZZ^c (n)), \QQ)
    \to \cdots
  \end{multline}
  Our calculations of motivic cohomology in
  Proposition~\ref{prop:structure-of-motivic-cohomology} give
  $$\Hom (H^i (X_\et, \ZZ^c (n)), \QQ) = 0 \quad \text{for }i \ne -1,$$
  and further by the definition of $\ZZ(n)$ in \eqref{eqn:definition-of-Z(n)},
  $$H^i_c (X_\et, \ZZ(n)) = 0 \quad \text{for }i \le 0.$$
  This implies that $H^i_\fg (X, \ZZ(n)) = 0$ for $i \le 0$.
  Since $H^i_c (G_\RR, X(\CC), \ZZ(n)) = 0$ for $i < 0$, we see from the exact
  sequence
  \begin{multline}
    \label{eqn:Wc-les}
    \cdots \to H^i_\Wc (X, \ZZ(n)) \to
    H^i_\fg (X,\ZZ(n)) \xrightarrow{H^i (i^*_\infty)}
    H^i_c (G_\RR, X(\CC), \ZZ(n)) \\
    \to H^{i+1}_\Wc (X, \ZZ(n)) \to \cdots
  \end{multline}
  that $H^i_\Wc (X, \ZZ(n)) = 0$ for $i \le 0$.

  \vspace{1em}

  For $i = 1$, the exact sequence \eqref{eqn:fg-les} shows that
  $H^1_c (X_\et, \ZZ(n)) \to H^1_\fg (X_\et, \ZZ(n))$ is an isomorphism.
  Consequently, we see that
  $\ker H^1 (i^*_\infty) \cong \ker H^1 (u_\infty^*)$:
  \[ \begin{tikzcd}
      H^1_c (X_\et, \ZZ(n)) \ar{r}{\cong}\ar{d}[swap]{H^1 (u^*_\infty)} & H^1_\fg (X, \ZZ(n)) \ar{dl}{H^1 (i^*_\infty)} \\
      H^1_c (G_\RR, X(\CC), \ZZ(n))
    \end{tikzcd} \]
  From long exact sequences \eqref{eqn:Wc-les} and
  \eqref{eqn:les-Hc-vs-H-hat-c}, we obtain short exact sequences
  \[ \begin{tikzcd}[column sep=1em,row sep=0pt]
      0 \ar{r} & H^0_c (G_\RR, X(\CC), \ZZ(n)) \ar{r} & H^1_\Wc (X, \ZZ(n)) \ar{r} & \ker H^1 (i^*_\infty) \ar{r} & 0 \\
      0 \ar{r} & \widehat{H}^0_c (G_\RR, X(\CC), \ZZ(n)) \ar{r} & \widehat{H}^1_c (X_\et, \ZZ(n)) \ar{r} & \ker H^1 (u_\infty^*) \ar{r} & 0
    \end{tikzcd} \]
  Since $\ker H^1 (i^*_\infty) \cong \ker H^1 (u_\infty^*)$,
  this is part 1) of the proposition.

  \vspace{1em}

  We proceed to compute $H^i_\Wc (X, \ZZ(n))$ for $i \ge 2$. It is more
  convenient to do this without passing explicitly through
  $H^i_\fg (X, \ZZ(n))$. Consider the morphism of complexes
  \[ \widehat{\alpha}_{X,n}\colon
    \RHom (R\Gamma (X_\et, \ZZ^c(n)), \QQ[-2]) \to
    R\widehat{\Gamma}_c (X_\et, \ZZ(n)), \]
  defined in the same way as $\alpha_{X,n}$ in \eqref{eqn:alpha-X-n-cohomology},
  only without the final projection from $\widehat{H}^i_c$ to $H^i_c$:
  \begin{multline*}
    H^i (\widehat{\alpha}_{X,n})\colon \Hom (H^{2-i} (X_\et, \ZZ^c(n)), \QQ)
    \xrightarrow{\QQ \twoheadrightarrow \QQ/\ZZ}
    H^{2-i} (X_\et, \ZZ^c(n))^D \\
    \xleftarrow{\cong} \widehat{H}^i_c (X_\et, \ZZ(n)).
  \end{multline*}
  The relation between $\widehat{\alpha}_{X,n}$ and $\alpha_{X,n}$
  is given by
  \[ \begin{tikzcd}
      \RHom (R\Gamma (X_\et, \ZZ^c(n)), \QQ[-2]) \ar{r}{\widehat{\alpha}_{X,n}}\ar{dr}[swap]{\alpha_{X,n}} & R\widehat{\Gamma}_c (X_\et, \ZZ(n)) \ar{d} \\
      & R\Gamma_c (X_\et, \ZZ(n))
    \end{tikzcd} \]
  Here the vertical arrow comes from the definition of modified \'{e}tale
  cohomology with compact support and it sits in an exact triangle
  \[ R\widehat{\Gamma}_c (X_\et, \ZZ(n)) \to
    R\Gamma_c (X_\et, \ZZ(n)) \xrightarrow{\widehat{u}^*_\infty}
    R\widehat{\Gamma}_c (G_\RR, X(\CC), \ZZ(n)) \to \cdots [1] \]
  ---see \cite[Lemma~6.14]{Flach-Morin-2018}.
  From the definition of $\widehat{\alpha}_{X,n}$ and the exact sequence
  \eqref{eqn:exact-sequence-A*-etc}, we calculate
  \begin{align*}
    \ker H^i (\widehat{\alpha}_{X,n}) & = H^{2-i} (X_\et, \ZZ^c(n))^*, \\
    \coker H^i (\widehat{\alpha}_{X,n}) & \cong (H^{2-i} (X_\et, \ZZ^c (n))_\tors)^D.
  \end{align*}

  We denote a cone of $\widehat{\alpha}_{X,n}$ by
  $R\widehat{\Gamma}_\fg (X, \ZZ(n))$ and set
  $$\widehat{H}^i_\fg (X, \ZZ(n)) \dfn H^i (R\widehat{\Gamma}_\fg (X, \ZZ(n))),$$
  so that there is a long exact sequence
  \begin{multline*}
    \cdots \to \Hom (H^{2-i} (X_\et, \ZZ^c (n)), \QQ) \xrightarrow{H^i (\widehat{\alpha}_{X,n})} \\
    \widehat{H}^i_c (X_\et, \ZZ(n)) \to
    \widehat{H}^i_\fg (X, \ZZ(n)) \to
    \Hom (H^{1-i} (X_\et, \ZZ^c (n)), \QQ) \to \cdots
  \end{multline*}
  The corresponding short exact sequences
  \[ 0 \to \coker H^i (\widehat{\alpha}_{X,n}) \to
    \widehat{H}^i_\fg (X, \ZZ(n)) \to
    \ker H^{i+1} (\widehat{\alpha}_{X,n}) \to 0 \]
  are split, since $\ker H^{i+1} (\widehat{\alpha}_{X,n})$ is a free
  group. Therefore, we have
  \[ \widehat{H}^i_\fg (X, \ZZ(n)) \cong
    H^{1-i} (X_\et, \ZZ^c(n))^*
    \oplus
    (H^{2-i} (X_\et, \ZZ^c (n))_\tors)^D. \]

  There is a commutative diagram with distinguished rows and columns
  \[ \begin{tikzcd}[font=\footnotesize,column sep=1em]
      \RHom (R\Gamma (X_\et, \ZZ^c(n)), \QQ[-2]) \ar{r}{\widehat{\alpha}_{X,n}} \ar{d}{id} &[0pt] R\widehat{\Gamma}_c (X_\et, \ZZ(n)) \ar{d} \ar{r} & R\widehat{\Gamma}_\fg (X, \ZZ(n)) \ar{r}\ar{d} &[-0.5em] {[+1]} \ar{d}{id} \\
      \RHom (R\Gamma (X_\et, \ZZ^c(n)), \QQ[-2]) \ar{r}{\alpha_{X,n}} \ar{d} & R\Gamma_c (X_\et, \ZZ(n)) \ar{r} \ar{d}{\widehat{u}^*_\infty} & R\Gamma_\fg (X, \ZZ(n)) \ar{r} \ar{d}{\widehat{i}^*_\infty} & {[+1]} \ar{d} \\
      0 \ar{r}\ar{d} & R\widehat{\Gamma}_c (G_\RR, X(\CC), \ZZ(n)) \ar{r}{id} \ar{d} & R\widehat{\Gamma}_c (G_\RR, X(\CC), \ZZ(n)) \ar{r} \ar{d} & 0 \ar{d} \\
      \RHom (R\Gamma (X_\et, \ZZ^c(n)), \QQ[-1]) \ar{r}{\widehat{\alpha}_{X,n} [1]} & R\widehat{\Gamma}_c (X_\et, \ZZ(n)) [1] \ar{r} & R\widehat{\Gamma}_\fg (X, \ZZ(n)) [1] \ar{r} & {[+2]}
    \end{tikzcd} \]
  Here $\widehat{u}^*_\infty$ (resp. $\widehat{i}^*_\infty$)
  is defined as the composition of the canonical morphism $u^*_\infty$
  (resp. $i^*_\infty$) with the projection to the Tate cohomology
  \[
    \pi\colon R\Gamma_c (G_\RR, X(\CC), \ZZ(n)) \to
    R\widehat{\Gamma}_c (G_\RR, X(\CC), \ZZ(n)).
  \]
  In our case of one-dimensional $X$, we know that $H^i (\pi)$ is an isomorphism
  for $i \ge 1$ (cf. \cite[Proposition~3.2]{Beshenov-Weil-etale-1}). Therefore,
  the five-lemma applied to
  \[ \begin{tikzcd}[font=\footnotesize]
      R\Gamma_\Wc (X, \ZZ(n)) \ar{r}\ar{d}{f} & R\Gamma_\fg (X, \ZZ(n)) \ar{r}{i^*_\infty}\ar{d}{id} & R\Gamma_c (G_\RR, X(\CC), \ZZ(n)) \ar{r}\ar{d}{\pi} & \cdots [1]\ar{d}{f [1]} \\
      R\widehat{\Gamma}_\fg (X, \ZZ(n)) \ar{r} & R\Gamma_\fg (X, \ZZ(n)) \ar{r}{\widehat{i}^*_\infty} & R\widehat{\Gamma}_c (G_\RR, X(\CC), \ZZ(n)) \ar{r} & \cdots [1]
    \end{tikzcd} \]
  shows that for $i \ge 2$ holds
  \[ H^i_\Wc (X, \ZZ(n)) \cong \widehat{H}^i_\fg (X, \ZZ(n))
    \cong H^{1-i} (X_\et, \ZZ^c(n))^*
    \oplus
    (H^{2-i} (X_\et, \ZZ^c (n))_\tors)^D. \]
  Our calculations of motivic cohomology in
  Proposition~\ref{prop:structure-of-motivic-cohomology} yield
  \begin{align*}
    H^2_\Wc (X, \ZZ(n)) & \cong H^{-1} (X_\et, \ZZ^c(n))^*
                          \oplus
                          H^0 (X_\et, \ZZ^c (n))^D, \\
    H^3_\Wc (X, \ZZ(n)) & \cong (H^{-1} (X_\et, \ZZ^c (n))_\tors)^D, \\
    H^i_\Wc (X, \ZZ(n)) & = 0 \text{ for }i \ge 4. \qedhere
  \end{align*}
\end{proof}

\begin{remark}
  A priori, the short exact sequence \eqref{eqn:H1Wc-ses} need not split.
  This will not bother us for the determinant calculations in
  \S\ref{sec:Weil-etale-proof} below.
\end{remark}

\begin{remark}
  The groups $H^i_\Wc (X, \ZZ(n))$ for $X = \Spec \mathcal{O}_F$ are already
  calculated in \cite[\S 5.8.3]{Flach-Morin-2018}. The result is (using the
  identification \eqref{eqn:Zc(n)-vs-Z(d-n)})
  \begin{equation}
    \label{eqn:Weil-etale-cohomology-of-Spec-OF}
    H^i_\Wc (X, \ZZ(n)) \cong
    \begin{cases}
      \ZZ^{\oplus d_n}, & i = 1, \\
      H^{-1} (X_\et, \ZZ^c (n))^* \oplus H^0 (X_\et, \ZZ^c (n))^D, & i = 2, \\
      (H^{-1} (X_\et, \ZZ^c (n))_\tors)^D, & i = 3, \\
      0, & i \ne 1,2,3.
    \end{cases}
  \end{equation}
  Our calculation generalizes this. What may look puzzling is the answer for
  $H^1_\Wc (X,\ZZ(n))$ given by Proposition~\ref{prop:calculation-of-H-Wc}.
  In the case of $X = \Spec \mathcal{O}_F$ we have, according to
  \eqref{eqn:2-torsion-in-H1-Zc-for-Spec-OF}, that
  $H^1 (X_\et, \ZZ^c(n)) \cong (\ZZ/2\ZZ)^{\oplus r_1}$ for even $n$, and
  hence $T_1 = 0$, which agrees with
  \eqref{eqn:Weil-etale-cohomology-of-Spec-OF}.

  Intuitively, the arithmetically interesting cohomology $H^i (X_\et, \ZZ^c(n))$
  for $X = \Spec \mathcal{O}_F$ is concentrated in degrees $i = -1,0$. The
  groups $H^i (X_\et, \ZZ^c (n))$ for $i \ge 1$ do not contain any interesting
  information: they are finite $2$-torsion, coming from the real places of
  $F$. The transition to Weil-\'{e}tale cohomology eliminates this
  $2$-torsion. On the other hand, the group $H^1 (X_\et, \ZZ^c (n))$ for a curve
  over a finite field $X/\FF_q$ is nontrivial and contains arithmetic
  information. The finite group $T_1$ appearing in the statement removes the
  $2$-torsion coming from the real places of $X$.
\end{remark}

\begin{remark}
  For a curve over a finite field $X/\FF_q$, all groups $H^i (X_\et, \ZZ^c(n))$
  are finite, and our calculation gives
  $H^i_\Wc (X, \ZZ(n)) \cong H^{2-i} (X_\et, \ZZ^c(n))^D$.
  This is true for any variety over a finite field $X/\FF_q$ and $n < 0$, under
  the assumption of finite generation of $H^i (X_\et, \ZZ^c(n))$; see
  \cite[Proposition~7.7]{Beshenov-Weil-etale-1}.
\end{remark}

\begin{remark}
  It is conjectured in \cite[\S 3]{Beshenov-Weil-etale-2} that
  \[ \ord_{s = n} \zeta (X,s) =
    \sum_{i\in \ZZ} (-1)^i\cdot i\cdot\rk_\ZZ H^i_\Wc (X,\ZZ(n)). \]
  In this case
  \begin{align*}
    \rk_\ZZ H^1_\Wc (X,\ZZ(n)) & = \rk_\ZZ H^2_\Wc (X,\ZZ(n)) = d_n, \\
    \rk_\ZZ H^3_\Wc (X,\ZZ(n)) & = 0,
  \end{align*}
  so the conjecture holds by Proposition~\ref{prop:vanishing-order-equals-dn}.
\end{remark}


\section{Weil-\'{e}tale proof of the special value formula}
\label{sec:Weil-etale-proof}

Now we explicitly write down the special value conjecture $\mathbf{C} (X,n)$
from \cite[\S 4]{Beshenov-Weil-etale-2}. To do this, consider the canonical
isomorphism
\begin{multline*}
  \lambda\colon \RR \xrightarrow{\cong}
  \bigotimes_{i\in \ZZ} (\det_\RR H^i_\Wc (X, \RR (n)))^{(-1)^i} \\
  \xrightarrow{\cong} \Bigl(\bigotimes_{i\in \ZZ} (\det_\ZZ H^i_\Wc (X, \ZZ (n)))^{(-1)^i}\Bigr) \otimes_\ZZ \RR \\
  \xrightarrow{\cong} (\det_\ZZ R\Gamma_\Wc (X, \ZZ (n))) \otimes_\ZZ \RR,
\end{multline*}
where the first isomorphism
$\RR \cong \bigotimes_{i\in \ZZ} (\det_\RR H^i_\Wc (X, \RR (n)))^{(-1)^i}$
comes from the regulator, as explained below.

\vspace{1em}

In our case, we are interested in the determinant of the Weil-\'{e}tale complex
\begin{multline*}
  \det_\ZZ R\Gamma_\Wc (X, \ZZ(n)) \cong
  \bigotimes_{i \in \ZZ} \det_\ZZ H^i_\Wc (X, \ZZ(n))^{(-1)^i} \\
  =
  \det_\ZZ H^1_\Wc (X, \ZZ(n))^{-1} \otimes_\ZZ
  \det_\ZZ H^2_\Wc (X, \ZZ(n)) \otimes_\ZZ
  \det_\ZZ H^3_\Wc (X, \ZZ(n))^{-1}.
\end{multline*}
Using the calculations from Proposition~\ref{prop:calculation-of-H-Wc},
\begin{align*}
  \det_\ZZ H^1_\Wc (X, \ZZ(n)) & \cong \det_\ZZ H^0_c (G_\RR, X(\CC), \ZZ(n)) \otimes_\ZZ \det_\ZZ T_1, \\
  \det_\ZZ H^2_\Wc (X, \ZZ(n)) & \cong \det_\ZZ H^{-1} (X_\et, \ZZ^c(n))^* \otimes_\ZZ \det_\ZZ H^0 (X_\et, \ZZ^c(n))^D, \\
  \det_\ZZ H^3_\Wc (X, \ZZ(n)) & \cong \det_\ZZ (H^{-1} (X_\et, \ZZ^c(n))_\tors)^D.
\end{align*}
So we have an isomorphism (up to sign $\pm 1$, after rearranging the terms)
\begin{multline*}
  \det_\ZZ R\Gamma_\Wc (X, \ZZ(n)) \cong \\
  =
  \det_\ZZ H^0_c (G_\RR, X(\CC), \ZZ(n))^{-1} \otimes_\ZZ
  \det_\ZZ H^{-1} (X_\et, \ZZ^c(n))^* \otimes_\ZZ \\
  \det_\ZZ (T_1)^{-1} \otimes_\ZZ  \det_\ZZ H^0 (X_\et, \ZZ^c(n))^D \otimes_\ZZ \det_\ZZ ((H^{-1} (X_\et, \ZZ^c(n))_\tors)^D)^{-1}.
\end{multline*}
Recall that $T_1$, $H^0 (X_\et, \ZZ^c(n))^D$,
$(H^{-1} (X_\et, \ZZ^c(n))_\tors)^D$ are finite groups, while the groups
$H^0_c (G_\RR, X(\CC), \ZZ(n))$ and $H^{-1} (X_\et, \ZZ^c(n))^*$ are free of
rank $d_n$. Now we consider the canonical trivialization
\[ (\det_\ZZ R\Gamma_\Wc (X, \ZZ(n))) \otimes_\ZZ \RR \cong
  \bigotimes_{i\in \ZZ} \det_\RR (H^i_\Wc (X, \ZZ(n))\otimes_\ZZ \RR)
  \cong \RR \]
via the regulator morphism
\[ \begin{tikzcd}[column sep=3em]
    H^0_c (G_\RR, X(\CC), \ZZ(n)) \otimes \RR\ar{d}{\cong} & \Hom (H^{-1} (X_\et, \ZZ^c (n)), \ZZ) \otimes \RR\ar{d}{\cong} \\
    H^0_c (G_\RR, X(\CC), \RR (n)) \ar{r}{Reg^\vee_{X,n}}[swap]{\cong} & \Hom (H^{-1} (X_\et, \ZZ^c (n)), \RR)
  \end{tikzcd} \]

\begin{proposition}
  Under the above trivialization, $\det_\ZZ R\Gamma_\Wc (X, \ZZ(n)) \subset \RR$
  corresponds to $\alpha^{-1}\,\ZZ \subset \RR$, where
  \begin{align*}
    \alpha & = \frac{|H^0 (X_\et, \ZZ^c(n))^D|}{|T_1| \cdot |(H^{-1} (X_\et, \ZZ^c(n))_\tors)^D|}\,R_{X,n} \\
           & = 2^\delta\,\frac{|H^0 (X_\et, \ZZ^c (n))|}{|H^{-1} (X_\et, \ZZ^c(n))_\tors|\cdot |H^1 (X_\et, \ZZ^c(n))|}\,R_{X,n},
  \end{align*}
  the number $\delta$ is given by \eqref{eqn:delta}, and $R_{X,n}$ is
  the regulator from Definition~\ref{dfn:regulator}.

  \begin{proof}
    For the finite groups $T_1$, $H^0 (X_\et, \ZZ^c(n))^D$,
    $(H^{-1} (X_\et, \ZZ^c(n))_\tors)^D$, this is
    \cite[Lemma~A.5]{Beshenov-Weil-etale-2}. For the free groups
    $H^0_c (G_\RR, X(\CC), \ZZ(n))$ and $H^{-1} (X_\et, \ZZ^c (n))^*$, on the
    other hand, this is Proposition~\ref{prop:trivialization-of-free-part}
    (now our groups sit in degrees $1$ and $2$, so the determinant gets
    inverted).
  \end{proof}
\end{proposition}

We recall that Conjecture $\mathbf{C} (X,n)$ from
\cite[\S 4]{Beshenov-Weil-etale-2} states that the canonical embedding
$\det_\ZZ R\Gamma_\Wc (X, \ZZ(n)) \subset \RR$ corresponds to
$\zeta^* (X,n)^{-1}\,\ZZ \subset \RR$.

\begin{proposition}
  Let $X$ be a one-dimensional arithmetic scheme and ${n < 0}$. Then the special
  value conjecture $\mathbf{C} (X,n)$ stated in \cite{Beshenov-Weil-etale-2} is
  equivalent to formula \eqref{eqn:special-value-formula}.
\end{proposition}

In \cite[\S 7]{Beshenov-Weil-etale-2} it is already proved (using essentially
the same localization idea as in this text) that $\mathbf{C} (X,n)$ holds
unconditionally for an abelian one-dimensional arithmetic scheme $X$. Together
with the proposition above, this proves Theorem~\ref{main-theorem} from the
introduction.


\section{A couple of examples}
\label{sec:examples}

We conclude with two examples that illustrate how localization arguments
work. The first is rather general and consists in specifying
\S\ref{sec:direct-proof} to the case of a non-maximal order in a number field.

\begin{example}
  Let $\mathcal{O} \subset \mathcal{O}_F$ be a non-maximal order in a number
  field $F/\QQ$. Denote $X = \Spec \mathcal{O}$ and
  $X' = \Spec \mathcal{O}_F$. Geometrically, $\nu\colon X' \to X$ is the
  normalization. There exist open dense subschemes $U \subset X$ and
  $U' \subset X'$ such that $\nu$ induces an isomorphism $U' \cong U$. If we
  denote the corresponding closed complements by $Z = X\setminus U$ and
  $Z' = X'\setminus U'$, then we have
  $$\zeta_\mathcal{O} (s) = \frac{\zeta (Z,s)}{\zeta (Z',s)}\,\zeta_F (s).$$
  For this identity formulated in classical terms of algebraic number theory,
  see, for example, \cite{Jenner-1969}. In particular,
  \[ \zeta^*_\mathcal{O} (n) =
    \pm \frac{|H^1 (Z'_\et, \ZZ^c(n))|}{|H^1 (Z_\et, \ZZ^c(n))|}\,\zeta^*_F (n). \]

  Now our special value conjectures for $\zeta^*_\mathcal{O} (n)$ and
  $\zeta^*_F (n)$ take the form
  \begin{align}
    \label{eqn:special-value-nonmax-order}
    \zeta^*_\mathcal{O} (n) & \stackrel{?}{=}
                              \pm 2^\delta\,\frac{|H^0 (X_\et, \ZZ^c (n))|}{|H^{-1} (X_\et, \ZZ^c (n))_\tors|\cdot |H^1 (X_\et, \ZZ^c (n))|}\,R, \\
    \label{eqn:special-value-max-order}
    \zeta^*_F (n) & \stackrel{?}{=}
                    \pm 2^\delta\,\frac{|H^0 (X'_\et, \ZZ^c (n))|}{|H^{-1} (X'_\et, \ZZ^c (n))_\tors|\cdot |H^1 (X'_\et, \ZZ^c (n))|}\,R.
  \end{align}
  Here
  $|H^{-1} (X_\et, \ZZ^c (n))_\tors| = |H^{-1} (X'_\et, \ZZ^c (n))_\tors|$, and
  the exact sequences of finite groups
  \begin{multline*}
    0 \to H^0 (X_\et, \ZZ^c (n)) \to H^0 (U_\et, \ZZ^c (n)) \to \\
    H^1 (Z_\et, \ZZ^c (n)) \to H^1 (X_\et, \ZZ^c (n)) \to H^1 (U_\et, \ZZ^c (n)) \to 0
  \end{multline*}
  \begin{multline*}
    0 \to H^0 (X'_\et, \ZZ^c (n)) \to H^0 (U'_\et, \ZZ^c (n)) \to \\
    H^1 (Z'_\et, \ZZ^c (n)) \to H^1 (X'_\et, \ZZ^c (n)) \to H^1 (U'_\et, \ZZ^c (n)) \to 0
  \end{multline*}
  give us
  \[ \frac{|H^1 (Z'_\et, \ZZ^c (n))|}{|H^1 (Z_\et, \ZZ^c (n))|} =
    \frac{|H^1 (X'_\et, \ZZ^c (n))|}{|H^1 (X_\et, \ZZ^c (n))|}\cdot
    \frac{|H^0 (X_\et, \ZZ^c (n))|}{|H^0 (X'_\et, \ZZ^c (n))|}, \]
  which implies that the formulas \eqref{eqn:special-value-nonmax-order} and
  \eqref{eqn:special-value-max-order} are equivalent.
\end{example}

The second example is suggested by \cite[\S 7]{Jordan-Poonen-2020}.

\begin{example}
  Let $p$ be an odd prime. Consider the affine scheme
  \[ X = \Spec (\ZZ [1/2] \times_{\FF_p} \FF_p [t]) =
    \Spec \ZZ [1/2] \mathop{\sqcup}_{\Spec \FF_p} \AA^1_{\FF_p} \]
  obtained from $\Spec \ZZ [1/p]$ and $\AA^1_{\FF_p} = \Spec \FF_p [t]$ by
  gluing together the points corresponding to the prime ideals
  $(p) \subset \ZZ [1/2]$ and $(t) \subset \FF_p [t]$:
  \[ \ZZ [1/2] \times_{\FF_p} \FF_p [t] =
    \{ (a,f) \in \ZZ [1/2] \times \FF_p [t] \mid a \equiv f (0) \pmod{p} \}. \]
  If we take odd $n < 0$, then there is no regulator. Let us consider $n = -3$.

  First, recall some calculations of the motivic cohomology of $\Spec \ZZ$.
  Using \cite[Proposition~2.1]{Kolster-Sands-2008} and known calculations of the
  $K$-groups of $\ZZ$ (see Weibel's survey \cite{Weibel-2005}), we get
  \begin{align*}
    H^{-1} (\Spec \ZZ_\et, \ZZ^c (-3)) & \cong K_7 (\ZZ) \cong \ZZ/240\ZZ, \\
    H^0 (\Spec \ZZ_\et, \ZZ^c (-3)) & \cong \ZZ/2\ZZ, \\
    H^1 (\Spec \ZZ_\et, \ZZ^c (-3)) & = 0.
  \end{align*}
  We note that, as expected,
  \[ \zeta (\Spec \ZZ, -3) = \zeta (-3) = -\frac{B_4}{4} = \frac{1}{120} =
    \frac{|H^0 (\Spec \ZZ_\et, \ZZ^c (-3))|}{|H^{-1} (\Spec \ZZ_\et, \ZZ^c (-3))|}. \]

  The localization gives
  \begin{align*}
    H^{-1} (\Spec \ZZ [1/2]_\et, \ZZ^c (-3)) & \cong H^{-1} (\Spec \ZZ_\et, \ZZ^c (-3)) \cong \ZZ/240\ZZ, \\
    H^0 (\Spec \ZZ [1/2]_\et, \ZZ^c (-3)) & \cong \ZZ/2\ZZ \oplus \ZZ/7\ZZ, \\
    H^1 (\Spec \ZZ [1/2]_\et, \ZZ^c (-3)) & = H^1 (\Spec \ZZ_\et, \ZZ^c (-3)) = 0.
  \end{align*}
  Arithmetically, this corresponds to the fact that the zeta function of
  $\Spec \ZZ [1/2]$ has the same Euler product as $\zeta (s)$, with the factor
  $\frac{1}{1-2^{-s}}$ removed. Therefore, when $s = -3$, the zeta-value should
  be corrected by $2^3 - 1 = 7$.

  For $\AA^1_{\FF_p}$, we now have
  \[ H^i (\AA^1_{\FF_p, \et}, \ZZ^c (n)) \cong
    H^{i+2} (\Spec \FF_{p,\et}, \ZZ^c (n-1)) \cong
    \begin{cases}
      \ZZ/(p^{1-n} - 1)\ZZ, & i = -1, \\
      0, & i \ne -1.
    \end{cases} \]
  In particular, the motivic cohomology of $\AA^1_{\FF_p}$ is concentrated in
  \[ H^{-1} (\AA^1_{\FF_p, \et}, \ZZ^c (-3)) \cong \ZZ/(p^4-1)\ZZ. \]

  Consider the normalization of $X$, given by
  $X' = \Spec \ZZ [1/2] \sqcup \AA^1_{\FF_p}$:
  \[ \begin{tikzcd}
      Z' \ar[right hook->]{r}\ar{d} & X' \ar{d} \\
      Z \ar[right hook->]{r} & X
    \end{tikzcd} \]
  Here $Z = \{ \mathfrak{p} \}$, $Z' = \{ \mathfrak{P}, \mathfrak{P}' \}$,
  and
  \begin{align*}
    \mathfrak{p} & \dfn \{ (a,f) \in \ZZ [1/2] \times \FF_p [t] \mid a \equiv f(0) \equiv 0 \pmod{p} \}, \\
    \mathfrak{P} & \dfn \{ (a,f) \in \ZZ [1/2] \times \FF_p [t] \mid a \equiv 0 \pmod{p} \}, \\
    \mathfrak{P}' & \dfn \{ (a,f) \in \ZZ [1/2] \times \FF_p [t] \mid f(0) \equiv 0 \pmod{p} \}.
  \end{align*}
  The canonical morphism $X' \to X$ induces an isomorphism
  \[ X'\setminus Z' \cong X\setminus Z \cong
    (\Spec \ZZ \setminus \{ (2), (p) \}) \sqcup
    (\Spec \FF_p [t] \setminus (t)). \]

  We calculate via localizations that
  \begin{align*}
    H^{-1} (X_\et, \ZZ^c (-3)) & \cong H^{-1} ((X\setminus Z)_\et, \ZZ^c (-3)) \cong
                                 \ZZ/240\ZZ \oplus \ZZ/(p^4 - 1)\ZZ, \\
    H^0 (X_\et, \ZZ^c (-3)) & \cong \ZZ/2\ZZ \oplus \ZZ/7\ZZ \oplus \ZZ/(p^3 - 1)\ZZ, \\
    H^1 (X_\et, \ZZ^c (-3)) & = 0.
  \end{align*}
  Consequently,
  \[ \frac{|H^0 (X_\et, \ZZ^c (-3))|}{|H^{-1} (X_\et, \ZZ^c (-3))|\cdot |H^1 (X_\et, \ZZ^c (-3))|} =
    \frac{7}{120}\,\frac{p^3-1}{p^4 - 1}. \]

  At the level of zeta-functions,
  \begin{multline*}
    \zeta (X,s) =
    \zeta (Z,s)\,\zeta (X\setminus Z,s) =
    \frac{\zeta (Z,s)}{\zeta (Z',s)}\,\zeta (X',s) \\
    = \frac{1}{\zeta (\Spec \FF_p,s)}\,\zeta (\Spec \ZZ[1/2],s)\,\zeta (\AA^1_{\FF_p},s) \\
    = (1-p^{-s})\,(1 - 2^{-s})\,\zeta (s)\,\frac{1}{1 - p^{1-s}}.
  \end{multline*}

  In particular, substituting $s = -3$, we get
  $$\zeta (X,-3) = -\frac{7}{120}\,\frac{p^3 - 1}{p^4 - 1}.$$
\end{example}


\pagebreak
\bibliographystyle{abbrv}
\bibliography{../weil-etale}

\end{document}